\documentclass[12pt, a4paper]{article}
\usepackage{rotating}
\usepackage{graphics}
\usepackage[all]{xy}
\xyoption{v2}
\xyoption{knot}

\author {Stefan Forcey}

\title {Vertically Iterated Classical Enrichment}

\newtheorem{theorem}{Theorem}
\newtheorem{definition}{Definition}

\newcommand{\MySection}[1]
{\section{ #1}}

\newcommand{\bcal}[1]{\mbox{\boldmath${\cal {#1}}$}}

\begin{document}
\SloppyCurves{

\maketitle
\begin{abstract}
 Lyubashenko has described enriched 2--categories as categories enriched over 
     ${\cal V}$--Cat, the 2--category of categories enriched over a symmetric monoidal ${\cal V}$.
This construction is the strict analogue for ${\cal V}$--functors in ${\cal V}$--Cat of Brian Day's probicategories
for ${\cal V}$--modules in ${\cal V}$--Mod.
 Here 
     I generalize the strict version to enriched $n$--categories for $k$--fold monoidal ${\cal V}$. The latter is defined 
     as by Balteanu, Fiedorowicz, Schw${\rm \ddot a}$nzl and Vogt but with the addition of making 
     visible the coherent associators  $\alpha^i$.
     The symmetric case can easily be recovered. This paper
     proposes a recursive definition of ${\cal V}$--$n$--categories and their morphisms.
     We show that 
      for ${\cal V}$ $k$--fold monoidal the structure of a $(k-n)$--fold monoidal strict 
    $(n+1)$--category is possessed by  ${\cal V}$--$n$--Cat. This article is a completion of the work begun 
    in \cite{forcey2}, and the initial sections duplicate the beginning of that paper.

\end{abstract}

      
            \MySection{Introduction}
            
	        There is an ongoing massive effort to link category theory and geometry, just a part 
	        of the broad undertaking known as categorification as described by Baez and Dolan in \cite{Baez1}.
	        This effort has as a partial goal that
	        of understanding
	        the categories and functors that correspond to loop spaces and their
	        associated topological functors. 
	        Progress towards this goal has been advanced greatly by the recent work of Balteanu,
	        Fiedorowicz, Schw${\rm \ddot a}$nzl, and Vogt in \cite{Balt}
	        where they show a direct correspondence between $k$-fold monoidal
	        categories and  $k$-fold loop spaces through the
	        categorical nerve. 
	        
	        As I pursued part of a plan to relate the enrichment functor to topology, I
	        noticed that the concept of higher dimensional
	        enrichment would be important in its relationship to double, triple and
	        further iterations of delooping.
	        The concept of enrichment over a monoidal category is well known, and 
	        enriching over the category of categories enriched over 
	        a monoidal category is  
	        defined, for the case of symmetric categories, in the paper on
	        $A_{\infty}$--categories
	        by Lyubashenko,  \cite{Lyub}. In the case of ${\cal V}$ closed, symmetric, and cocomplete 
              these are equivalent to the probicategories described in a preprint of Day. The latter are 
many-object versions of special cases of the
              promonoidal categories Day defines in \cite{Day}.
	    It seems that it is a good idea to generalize Lyubashenko's definition first
	        to the case of an iterated monoidal
	        base category
	        and then to 
	        define ${\cal V}$--$(n+1)$--categories as
	        categories enriched over ${\cal V}$--$n$--Cat, the $(k-n)$--fold monoidal
	        strict 
	        $(n+1)$--category of ${\cal V}$--$n$--categories where $k>n \in {\mathbf N}$. Of course the facts
	        implicit in this last statement must be verified. 

               The adjective ``vertical'' in the title is meant to distinguish this sort of iteration
from another possibility: that of restricting attention to the monoids in ${\cal V}$--Cat. These are 
monoidal ${\cal V}$--categories, or one-object ${\cal V}$--2--categories. Also these can be viewed as
enriched categories with monoidal underlying categories, and enrichment over a monoidal ${\cal V}$--category turns out
to be enrichment over that underlying category. Iteration of the enrichment in that sense will be arbitrarily referred to as
horizontal.
For now we consider the ``perpendicular'' direction.
 
	        At each stage of successive enrichments, the number of monoidal products
	        should decrease and the categorical dimension
	        should increase, both by one. This is motivated by topology. When we
	        consider the loop space of a topological space,
	        we see that paths (or  1--cells) in the original 
	        are now points (or objects) in the derived space. There is also now
	        automatically a product structure
	        on the points in the derived space, where multiplication is given by
	        concatenation of loops. Delooping is the inverse 
	        functor here, and thus involves shifting objects to the status of
	        1--cells and decreasing the number of ways to multiply.

	        The concept of a $k$--fold monoidal strict $n$--category  
	        is easy enough to define as a tensor object in a category of $(k-1)$--fold
	        monoidal $n$--categories with cartesian product. Thus the products and
	        accompanying associator and interchange transformations
	        are strict $n$--functors and $n$--natural transformations respectively. That
	        this sort of structure 
	        ($(k-n)$--fold monoidal strict 
	        $n+1$ category) is possessed by  ${\cal V}$--$n$--Cat for ${\cal V}$ $k$--fold monoidal is shown 
	        for $n = 1$ and all $k$ in my paper \cite{forcey1}.
              The case $n = 2$ is shown in \cite{forcey2} 
	        This paper completes the equation by presenting a
	        full inductive proof covering all $n,k.$ 
	        
	        In general the decrease is engineered by a shift in index--we define 
	        new products ${\cal V}$--$n$--Cat$\times {\cal V}$--$n$--Cat $\to {\cal V}$--$n$--Cat by 
	        using cartesian products of object sets and letting hom--objects of the
	        $i$th product of enriched $n$--categories
	        be the $(i+1)$th product of hom--objects of the component categories.
	        Symbolically,
	        $$({\cal A} \otimes^{(n)}_{i} {\cal B})((A,B),(A',B')) = {\cal A}(A,A')\otimes^{(n-1)}_{i+1} {\cal B}(B,B').$$
	        The superscript $(n)$ is not necessary since the product is defined by
	        context, 
	                but I insert it to make clear at what
	                level of enrichment the product is occurring. 
	        Defining the necessary natural transformations for this new product as
	        ``based upon'' the old ones, and the checking of 
	        the axioms that define their structure is briefly mentioned later on in
	        this paper and more fully  described in \cite{forcey1}  
	         for certain cases.
	          
	        The definition of a category enriched over ${\cal V}$--$n$--Cat is simply
	        stated by describing the process as enriching over 
    ${\cal V}$--$n$--Cat with the first of the $k-n$ ordered products.
        In section 2 we quickly review the necessary preliminary definitions just as in \cite{forcey2}.  In section 3 we define 
        ${\cal V}$--$n$--categories and ${\cal V}$--$n$--functors, and in section 4 
        we discuss the change of base of enrichment 
    in the $k$--fold monoidal context and the forgetful functors thus derived. In section 5 we apply these results to 
    prove a general theorem about the categorical dimension of ${\cal V}$--$n$--Cat and describe the 
    specific higher morphisms that exist.
    
    \clearpage
          \newpage
        \MySection{Review of Definitions}
        
     In this section I briefly review the definitions of a category enriched over a monoidal category ${\cal V}$, a
     category enriched over an iterated monoidal category, and an enriched 2--category.
     I begin with the basic definitions of enrichment, included due to how often they 
      are referred to and followed as models in the rest of the paper. This first set of definitions
      can be found with more detail in \cite{Kelly} and \cite{EK1}.
      
      \begin{definition} For our purposes a {\it monoidal category} is a category ${\cal V}$
        together with a functor 
        $\otimes: {\cal V}\times{\cal V}\to{\cal V}$  and an object $I$ such that
        \begin{enumerate}
        \item $\otimes$ is  associative up to the coherent natural transformations $\alpha$. The coherence
        axiom is given by the commuting pentagon
        
        \noindent
	      	\begin{center}
	      	\resizebox{4.5in}{!}{
	              $$
	              \xymatrix@R=45pt@C=-35pt{
	              &((U\otimes V)\otimes W)\otimes X \text{ }\text{ }
	              \ar[rr]^{ \alpha_{UVW}\otimes 1_{X}}
	              \ar[dl]|{ \alpha_{(U\otimes V)WX}}
	              &\text{ }\text{ }\text{ }\text{ }\text{ }\text{ }\text{ }\text{ }\text{ }&\text{ }\text{ }(U\otimes (V\otimes W))\otimes X
	              \ar[dr]|{ \alpha_{U(V\otimes W)X}}&\\
	              (U\otimes V)\otimes (W\otimes X)
	              \ar[drr]|{ \alpha_{UV(W\otimes X)}}
	              &&&&U\otimes ((V\otimes W)\otimes X)
	              \ar[dll]|{ 1_{U}\otimes \alpha_{VWX}}
	              \\&&U\otimes (V\otimes (W\otimes X))&&&
	             }
	             $$
	             }
		                    \end{center}
	                   	                    
			       \item $I$ is a strict $2$-sided unit for $\otimes$.
			       \end{enumerate}
			       \end{definition}
			       \begin {definition} \label{V:Cat} A (small) ${\cal V}$ {\it--Category} ${\cal A}$ is a set $\left|{\cal A}\right|$ of 
			       {\it objects}, 
			       a {\it hom-object} ${\cal A}(A,B) \in \left|{\cal V}\right|$ for
			       each pair of objects of ${\cal A}$, a family of {\it composition morphisms} $M_{ABC}:{\cal A}(B,C)
			       \otimes{\cal A}(A,B)\to{\cal A}(A,C)$ for each triple of objects, and an {\it identity element} $j_{A}:I\to{\cal A}(A,A)$ for each object.
			       The composition morphisms are subject to the associativity axiom which states that the following pentagon commutes

			             \noindent
			             	          \begin{center}
			       	          \resizebox{4.5in}{!}{
			             $$
			             \xymatrix@R=45pt@C=-22pt{
			             &({\cal A}(C,D)\otimes {\cal A}(B,C))\otimes {\cal A}(A,B)\text{ }\text{ }
			             \ar[rr]^{\scriptstyle \alpha}
			             \ar[dl]^{\scriptstyle M \otimes 1}
			             &\text{ }\text{ }\text{ }\text{ }\text{ }\text{ }\text{ }\text{ }\text{ }\text{ }\text{ }\text{ }\text{ }\text{ }\text{ }&\text{ }\text{ }{\cal A}(C,D)\otimes ({\cal A}(B,C)\otimes {\cal A}(A,B))
			             \ar[dr]^{\scriptstyle 1 \otimes M}&\\
			             {\cal A}(B,D)\otimes {\cal A}(A,B)
			             \ar[drr]^{\scriptstyle M}
			             &&&&{\cal A}(C,D)\otimes {\cal A}(A,C)
			             \ar[dll]^{\scriptstyle M}
			             \\&&{\cal A}(A,D))&&&
			             }$$
			             }
			            	                    \end{center}
			        	                    
			       and to the unit axioms which state that both the triangles in the following diagram commute
			       
			       $$
			         \xymatrix{
			         I\otimes {\cal A}(A,B)
			         \ar[rrd]^{=}
			         \ar[dd]_{j_{B}\otimes 1}
			         &&&&{\cal A}(A,B)\otimes I 
			         \ar[dd]^{1\otimes j_{A}}
			         \ar[lld]^{=}\\
			         &&{\cal A}(A,B)\\
			         {\cal A}(B,B)\otimes {\cal A}(A,B)
			         \ar[rru]^{M_{ABB}}
			         &&&&{\cal A}(A,B)\otimes {\cal A}(A,A)
			         \ar[llu]^{M_{AAB}}
			         }
			      $$
			      
			      \end{definition}
			      
			       In general a ${\cal V}$--category is directly analogous to an (ordinary) category enriched over $\mathbf{Set}.$  
			       If ${\cal V} = \mathbf{Set}$ then these diagrams are the usual category axioms.

			      \begin{definition} \label{enriched:funct} For ${\cal V}$--categories 
			      ${\cal A}$ and ${\cal B}$, a ${\cal V}$--$functor$ $T:{\cal A}\to{\cal B}$ is a function
			       $T:\left| {\cal A} \right| \to \left| {\cal B} \right|$ and a family of 
			       morphisms $T_{AB}:{\cal A}(A,B) \to {\cal B}(TA,TB)$ in ${\cal V}$ indexed by 
			       pairs $A,B \in \left| {\cal A} \right|$.
			       The usual rules for a functor that state $T(f \circ g) = Tf \circ Tg$ 
			       and $T1_{A} = 1_{TA}$ become in the enriched setting, respectively, the commuting diagrams
			      
			      $$
			       \xymatrix{
			       &{\cal A}(B,C)\otimes {\cal A}(A,B)
			       \ar[rr]^{\scriptstyle M}
			       \ar[d]^{\scriptstyle T \otimes T}
			       &&{\cal A}(A,C)
			       \ar[d]^{\scriptstyle T}&\\
			       &{\cal B}(TB,TC)\otimes {\cal B}(TA,TB)
			       \ar[rr]^{\scriptstyle M}
			       &&{\cal B}(TA,TC)
			       }
			      $$
			     and
			      $$
			       \xymatrix{
			       &&{\cal A}(A,A)
			       \ar[dd]^{\scriptstyle T_{AA}}\\
			       I
			       \ar[rru]^{\scriptstyle j_{A}}
			       \ar[rrd]_{\scriptstyle j_{TA}}\\
			       &&{\cal B}(TA,TA).
			       }
			      $$
			     ${\cal V}$--functors can be composed to form a category called ${\cal V}$--Cat. This category
			     is actually enriched over $\mathbf{Cat}$, the category of (small) categories with cartesian product. 
			        \end{definition}
			     \begin{definition} \label{enr:nat:trans}
			     For ${\cal V}$--functors $T,S:{\cal A}\to{\cal B}$ a ${\cal V}$--{\it natural  
			      transformation} $\alpha:T \to S:{\cal A} \to {\cal B}$
			     is an $\left| {\cal A} \right|$--indexed family of 
			     morphisms $\alpha_{A}:I \to {\cal B}(TA,SA)$ satisfying the ${\cal V}$--naturality
			     condition expressed by the commutativity of the following hexagonal diagram:
			     
			     \noindent
			     	            	          \begin{center}
			         	          \resizebox{4.5in}{!}{
			      $$
			       \xymatrix{
			       &I \otimes {\cal A}(A,B)
			       \ar[rr]^-{\scriptstyle \alpha_{B} \otimes T_{AB}}
			       &&{\cal B}(TB,SB) \otimes {\cal B}(TA,TB)
			       \ar[rd]^-{\scriptstyle M}
			     \\
			       {\cal A}(A,B)
			       \ar[ru]^{=}
			     \ar[rd]_{=}
			       &&&&{\cal B}(TA,SB)
			     \\
			       &{\cal A}(A,B) \otimes I
			       \ar[rr]_-{\scriptstyle S_{AB} \otimes \alpha_{A}}
			       &&{\cal B}(SA,SB) \otimes {\cal B}(TA,SA)
			       \ar[ru]^-{\scriptstyle M}
			       }
			      $$
			     }
			            	                    \end{center}
			      \end{definition}
			      
			      For two ${\cal V}$--functors  $T,S$ to be equal is to say $TA = SA$ for all $A$ 
			      and for the ${\cal V}$--natural isomorphism $\alpha$ between them to have 
			      components $\alpha_{A} = j_{TA}$. This latter implies equality of the hom--object morphisms: 
			      $T_{AB} = S_{AB}$ for all pairs of objects. The implication is seen by combining the second diagram in 
			      Definition \ref{V:Cat} with all the diagrams
			     in Definitions \ref{enriched:funct} and \ref{enr:nat:trans}.

    The fact that ${\cal V}$--Cat has the structure of a 2--category is demonstrated in \cite{Kelly}. Now
    we review the transfer to enriching over a $k$--fold monoidal category. The latter sort of category was
    developed and defined in \cite{Balt}. The authors describe its structure as arising from its description
    as a monoid in the
        category of $(k-1)$--fold monoidal categories. Here is that definition altered only slightly 
    to make visible the coherent associators as in \cite{forcey1}. In that paper I describe its structure as 
    arising from its description
    as a tensor object in the
        category of $(k-1)$--fold monoidal categories.
        
\begin{definition} An $n${\it -fold monoidal category} is a category ${\cal V}$
          with the following structure. 
          \begin{enumerate}
          \item There are $n$ distinct multiplications
          $$\otimes_1,\otimes_2,\dots, \otimes_n:{\cal V}\times{\cal V}\to{\cal V}$$
          for each of which the associativity pentagon commutes
          
          \noindent
  		          \begin{center}
  	          \resizebox{4.5in}{!}{
          $$
          \xymatrix@R=45pt@C=-35pt{
          &((U\otimes_i V)\otimes_i W)\otimes_i X \text{ }\text{ }
          \ar[rr]^{ \alpha^{i}_{UVW}\otimes_i 1_{X}}
          \ar[dl]^{ \alpha^{i}_{(U\otimes_i V)WX}}
          &\text{ }\text{ }\text{ }\text{ }\text{ }&\text{ }\text{ }(U\otimes_i (V\otimes_i W))\otimes_i X
          \ar[dr]^{ \alpha^{i}_{U(V\otimes_i W)X}}&\\
          (U\otimes_i V)\otimes_i (W\otimes_i X)
          \ar[drr]^{ \alpha^{i}_{UV(W\otimes_i X)}}
          &&&&U\otimes_i ((V\otimes_i W)\otimes_i X)
          \ar[dll]^{ 1_{U}\otimes_i \alpha^{i}_{VWX}}
          \\&&U\otimes_i (V\otimes_i (W\otimes_i X))&&&
          }
          $$
          }
  		                    \end{center}

          ${\cal V}$ has an object $I$ which is a strict unit
          for all the multiplications.
          \item For each pair $(i,j)$ such that $1\le i<j\le n$ there is a natural
          transformation
          $$\eta^{ij}_{ABCD}: (A\otimes_j B)\otimes_i(C\otimes_j D)\to
          (A\otimes_i C)\otimes_j(B\otimes_i D).$$
          \end{enumerate}
          These natural transformations $\eta^{ij}$ are subject to the following conditions:
          \begin{enumerate}
          \item[(a)] Internal unit condition: 
          $\eta^{ij}_{ABII}=\eta^{ij}_{IIAB}=1_{A\otimes_j B}$
          \item[(b)] External unit condition:
          $\eta^{ij}_{AIBI}=\eta^{ij}_{IAIB}=1_{A\otimes_i B}$
          \item[(c)] Internal associativity condition: The following diagram commutes
          
          \noindent
	  	            	          \begin{center}
	      	          \resizebox{5in}{!}{
	     $$
            \diagram
            ((U\otimes_j V)\otimes_i (W\otimes_j X))\otimes_i (Y\otimes_j Z)
            \xto[rrr]^{\eta^{ij}_{UVWX}\otimes_i 1_{Y\otimes_j Z}}
            \ar[d]^{\alpha^i}
            &&&\bigl((U\otimes_i W)\otimes_j(V\otimes_i X)\bigr)\otimes_i (Y\otimes_j Z)
            \dto^{\eta^{ij}_{(U\otimes_i W)(V\otimes_i X)YZ}}\\
            (U\otimes_j V)\otimes_i ((W\otimes_j X)\otimes_i (Y\otimes_j Z))
            \dto^{1_{U\otimes_j V}\otimes_i \eta^{ij}_{WXYZ}}
            &&&((U\otimes_i W)\otimes_i Y)\otimes_j((V\otimes_i X)\otimes_i Z)
            \ar[d]^{\alpha^i \otimes_j \alpha^i}
            \\
            (U\otimes_j V)\otimes_i \bigl((W\otimes_i Y)\otimes_j(X\otimes_i Z)\bigr)
            \xto[rrr]^{\eta^{ij}_{UV(W\otimes_i Y)(X\otimes_i Z)}}
            &&& (U\otimes_i (W\otimes_i Y))\otimes_j(V\otimes_i (X\otimes_i Z))
            \enddiagram
            $$
           }
	            	                    \end{center}
  
           \item[(d)] External associativity condition: The following diagram commutes
            
            \noindent
  	            	          \begin{center}
      	          \resizebox{5in}{!}{ 
             $$
             \diagram
            ((U\otimes_j V)\otimes_j W)\otimes_i ((X\otimes_j Y)\otimes_j Z)
            \xto[rrr]^{\eta^{ij}_{(U\otimes_j V)W(X\otimes_j Y)Z}}
            \ar[d]^{\alpha^j \otimes_i \alpha^j}
            &&& \bigl((U\otimes_j V)\otimes_i (X\otimes_j Y)\bigr)\otimes_j(W\otimes_i Z)
            \dto^{\eta^{ij}_{UVXY}\otimes_j 1_{W\otimes_i Z}}\\
            (U\otimes_j (V\otimes_j W))\otimes_i (X\otimes_j (Y\otimes_j Z))
            \dto^{\eta^{ij}_{U(V\otimes_j W)X(Y\otimes_j Z)}}
            &&&((U\otimes_i X)\otimes_j(V\otimes_i Y))\otimes_j(W\otimes_i Z)
            \ar[d]^{\alpha^j}
            \\
            (U\otimes_i X)\otimes_j\bigl((V\otimes_j W)\otimes_i (Y\otimes_j Z)\bigr)
            \xto[rrr]^{1_{U\otimes_i X}\otimes_j\eta^{ij}_{VWYZ}}
            &&& (U\otimes_i X)\otimes_j((V\otimes_i Y)\otimes_j(W\otimes_i Z))
            \enddiagram
            $$
          }
	           	                    \end{center}
   
          \item[(e)] Finally it is required  for each triple $(i,j,k)$ satisfying
          $1\le i<j<k\le n$ that
          the giant hexagonal interchange diagram commutes.
         \end{enumerate}
          
          \noindent
  		          \begin{center}
  	          \resizebox{6in}{!}{
          $$
          \xymatrix@C=-118pt{
          &((A\otimes_k A')\otimes_j (B\otimes_k B'))\otimes_i((C\otimes_k C')\otimes_j (D\otimes_k D'))
          \ar[ddl]|{\eta^{jk}_{AA'BB'}\otimes_i \eta^{jk}_{CC'DD'}}
          \ar[ddr]|{\eta^{ij}_{(A\otimes_k A')(B\otimes_k B')(C\otimes_k C')(D\otimes_k D')}}
          \\\\
          ((A\otimes_j B)\otimes_k (A'\otimes_j B'))\otimes_i((C\otimes_j D)\otimes_k (C'\otimes_j D'))
          \ar[dd]|{\eta^{ik}_{(A\otimes_j B)(A'\otimes_j B')(C\otimes_j D)(C'\otimes_j D')}}
          &&((A\otimes_k A')\otimes_i (C\otimes_k C'))\otimes_j((B\otimes_k B')\otimes_i (D\otimes_k D'))
          \ar[dd]|{\eta^{ik}_{AA'CC'}\otimes_j \eta^{ik}_{BB'DD'}}
          \\\\
          ((A\otimes_j B)\otimes_i (C\otimes_j D))\otimes_k((A'\otimes_j B')\otimes_i (C'\otimes_j D'))
          \ar[ddr]|{\eta^{ij}_{ABCD}\otimes_k \eta^{ij}_{A'B'C'D'}}
          &&((A\otimes_i C)\otimes_k (A'\otimes_i C'))\otimes_j((B\otimes_i D)\otimes_k (B'\otimes_i D'))
          \ar[ddl]|{\eta^{jk}_{(A\otimes_i C)(A'\otimes_i C')(B\otimes_i D)(B'\otimes_i D')}}
          \\\\
          &((A\otimes_i C)\otimes_j (B\otimes_i D))\otimes_k((A'\otimes_i C')\otimes_j (B'\otimes_i D'))
          }
          $$
          }
  		                    \end{center}

        \end{definition}

     The authors of \cite{Balt} remark that a symmetric monoidal category is $n$-fold monoidal for all $n$.
     This they demonstrate by letting
        $$\otimes_1=\otimes_2=\dots=\otimes_n=\otimes$$
        and defining (associators added by myself)
        $$\eta^{ij}_{ABCD}=\alpha^{-1}\circ (1_A\otimes \alpha)\circ (1_A\otimes (c_{BC}\otimes 1_D))\circ (1_A\otimes \alpha^{-1})\circ \alpha$$
        for all $i<j$. Here $c_{BC}: B\otimes C \to C\otimes B$ is the symmetry natural transformation.
     This provides the hint that enriching over a $k$--fold monoidal category is precisely a generalization of enriching over
     a symmetric category. In the symmetric case, to define a product in ${\cal V}$--Cat, we need $c_{BC}$ in order to create a
     middle exchange morphism $m$.  
     To describe products in ${\cal V}$--Cat for ${\cal V}$ $k$--fold monoidal we simply use $m=\eta$.
     
Before treating the general case of enriching over the $k$--fold monoidal category of enriched $n$--categories we examine the 
definition in the two lowest categorical dimensions. This will enlighten the following discussion. The careful
unfolding of the definitions here will stand in for any in depth unfolding of the same enriched constructions at higher levels.  
Categories enriched over $k$--fold monoidal ${\cal V}$ are carefully defined in \cite{forcey1}, where they  
     are shown to be the objects of a $(k-1)$--fold monoidal 2--category. Here we need only the definitions. Simply put,
     a category enriched over a $k$--fold monoidal ${\cal V}$ is a category enriched in the usual sense over $({\cal V}, \otimes_1, I, \alpha).$
     The other $k-1$ products in ${\cal V}$ are used up in the structure of ${\cal V}$--Cat. I will always denote the product(s) in ${\cal V}$--Cat
       with a superscript in parentheses that corresponds to the level of enrichment of the components of their domain.
       The product(s) in ${\cal V}$ 
    should logically then have a superscript (0) but I have suppressed this for brevity and to agree with my sources.
     For ${\cal V}$ $k$--fold monoidal we define the $i$th product of ${\cal V}$--categories 
     ${\cal A} \otimes^{(1)}_{i}{\cal B}$
        to have objects $\in \left|{\cal A}\right|\times \left|{\cal B}\right|$
        and to have hom--objects $\in \left|{\cal V}\right|$ given by
        
        $$({\cal A} \otimes^{(1)}_{i} {\cal B})((A,B),(A',B')) = {\cal A}(A,A')\otimes_{i+1} {\cal B}(B,B').$$
        
        Immediately we see that ${\cal V}$--Cat is $(k-1)$--fold monoidal by definition. (The full proof of
        this is in \cite{forcey1}.)
        The composition morphisms are 
        \begin{footnotesize}
	$$M_{(A,B)(A',B')(A'',B'')} : ({\cal A} \otimes^{(1)}_{i}{\cal
        B})((A',B'),(A'',B''))\otimes_{1}({\cal A} \otimes^{(1)}_{i}{\cal
        B})((A,B),(A',B'))\to ({\cal A} \otimes^{(1)}_{i}{\cal B})((A,B),(A'',B''))$$
        \end{footnotesize}
        given by
        $$
        \xymatrix{
        ({\cal A} \otimes^{(1)}_{i}{\cal B})((A',B'),(A'',B''))\otimes_{1}({\cal A} \otimes^{(1)}_{i}{\cal B})((A,B),(A',B'))
        \ar@{=}[d]\\
        ({\cal A}(A',A'')\otimes_{i+1}{\cal B}(B',B''))\otimes_{1}({\cal A}(A,A')\otimes_{i+1}{\cal B}(B,B'))
        \ar[d]_{\eta^{1,i+1}}\\
        ({\cal A}(A',A'')\otimes_{1}{\cal A}(A,A'))\otimes_{i+1}({\cal B}(B',B'')\otimes_{1}{\cal B}(B,B'))
        \ar[d]_{M_{AA'A''}\otimes_{2}M_{BB'B''}}\\
        ({\cal A}(A,A'')\otimes_{i+1}{\cal B}(B,B''))
        \ar@{=}[d]\\
        ({\cal A} \otimes^{(1)}_{i}{\cal B})((A,B),(A'',B''))
       }
    $$
        
        The identity element is given by $j_{(A,B)} = \xymatrix{I = I \otimes_{i+1} I
         							\ar[d]^{j_A \otimes_{i+1} j_B}
         							\\{\cal A}(A,A)\otimes_{i+1} {\cal B}(B,B)
         							\ar@{=}[d]
         							\\({\cal A} \otimes^{(1)}_{i}{\cal B})((A,B),(A,B))
     							}$
      
      The unit object in ${\cal V}$--1--Cat $={\cal V}$--Cat is the enriched category ${\bcal I}^{(1)}={\cal I}$ where $\left|{\cal I}\right| = \{0\}$ and 
       ${\cal I}(0,0) = I$. Of course $M_{000} = 1_I = j_0.$
       
       That each product $\otimes^{(1)}_i$ thus defined is a 2--functor  ${\cal V}$--Cat $\times$ ${\cal V}$--Cat $\to$ ${\cal V}$--Cat
       is seen easily. Its action on enriched functors and natural transformations is to form formal products using $\otimes_{i+1}$ of their
       associated morphisms. That the result of this action is a valid enriched functor or natural transformation always follows from the
       naturality of $\eta.$
       
       Associativity in ${\cal V}$--Cat
          must hold for each $\otimes^{(1)}_{i}$. The components of the 2--natural
          isomorphism $\alpha^{(1)i}$
          $$\alpha^{(1)i}_{{\cal A}{\cal B}{\cal C}}: ({\cal A} \otimes^{(1)}_{i} {\cal B})
           \otimes^{(1)}_{i} {\cal C} \to {\cal A} \otimes^{(1)}_{i} ({\cal B} \otimes^{(1)}_{i}
          {\cal C})$$
          are ${\cal V}$--functors
          that send ((A,B),C) to (A,(B,C)) and whose hom-components 
          \begin{footnotesize}              
          $$\alpha^{(1)i}_{{{\cal A}{\cal B}{\cal C}}_{((A,B),C)((A',B'),C')}}: [({\cal A} \otimes^{(1)}_{i} {\cal B}) \otimes^{(1)}_{i} {\cal C}](((A,B),C),((A',B'),C'))
          \to [{\cal A} \otimes^{(1)}_{i} ({\cal B} \otimes^{(1)}_{i} {\cal C})]((A,(B,C)),(A',(B',C')))$$
           \end{footnotesize}
   are given by 
          $$\alpha^{(1)i}_{{{\cal A}{\cal B}{\cal C}}_{((A,B),C)((A',B'),C')}}
          = \alpha^{i+1}_{{\cal A}(A,A'){\cal B}(B,B'){\cal C}(C,C')}.$$
               
     Now for the interchange 2--natural transformations $\eta^{(1)ij}$
       for $j\ge i+1$.
       We define the component morphisms $\eta^{(1)i,j}_{{\cal A}{\cal B}{\cal C}{\cal D}}$
       that make a 2--natural transformation between 2--functors. Each component must be an enriched functor.
       Their action on objects
       is to send $$((A,B),(C,D)) \in \left|({\cal A}\otimes^{(1)}_{j} {\cal B})\otimes^{(1)}_{i} ({\cal C}\otimes^{(1)}_{j} {\cal D})\right|$$
        to $$((A,C),(B,D)) \in \left|({\cal A}\otimes^{(1)}_{i} {\cal C})\otimes^{(1)}_{j} ({\cal B}\otimes^{(1)}_{i} {\cal D})\right|.$$
       The hom--object morphisms are given by
        $$\eta^{(1)i,j}_{{{\cal A}{\cal B}{\cal C}{\cal D}}_{(ABCD)(A'B'C'D')}} =
   \eta^{i+1,j+1}_{{\cal A}(A,A'){\cal B}(B,B'){\cal C}(C,C'){\cal D}(D,D')}.$$
       That the axioms regarding the associators and interchange transformations are all obeyed is established
       in \cite{forcey1}.
       
       We now describe categories enriched over 
     ${\cal V}$--Cat. These are defined for the symmetric case in \cite{Lyub}. Here 
    the definition of 
       ${\cal V}$--2--category is generalized for ${\cal V}$ a $k$--fold monoidal category with $k\ge 2.$ 
       The definition for symmetric monoidal ${\cal V}$ can be easily recovered just by letting $\otimes_1=\otimes_2=\otimes$,
       $\alpha^2=\alpha^1 = \alpha$ and $\eta=m.$ 
       
       {\bf Example } 
 A (small, strict) ${\cal V}$--{\it 2--category} ${\bcal U}$ consists of
       
       \begin{enumerate}
           \item A set of objects $\left|{\bcal U}\right|$
           \item For each pair of objects $A,B \in \left|{\bcal U}\right|$ a ${\cal V}$--category ${\bcal U}(A,B).$
       
       Of course then ${\bcal U}(A,B)$ consists of a set of objects (which play the role of the 1--cells in a 2--category) 
       and for each pair $f,g \in \left|{\bcal U}(A,B)\right|$ an object 
       ${\bcal U}(A,B)(f,g) \in {\cal V}$ (which plays the role 
       of the hom--set of 2--cells in a 2--category.) Thus the vertical composition morphisms of these $hom_{2}$--objects are in ${\cal V}$:
       $$M_{fgh}:{\bcal U}(A,B)(g,h) \otimes_{1} {\bcal U}(A,B)(f,g) \to 
       {\bcal U}(A,B)(f,h)$$
       
       Also, the vertical identity for a 1-cell object $a \in \left|{\bcal U}(A,B)\right|$ is  $j_{a} : I \to {\bcal U}(A,B)(a,a)$.
       The associativity and the units of vertical composition are then those given by the respective axioms of enriched categories.  
           \item For each triple of objects $A,B,C \in 
       \left|{\bcal U}\right|$ a ${\cal V}$--functor
       $${\cal M}_{ABC}:{\bcal U}(B,C) \otimes^{(1)}_{1} {\bcal U}(A,B) \to {\bcal U}(A,C)$$ 
       Often I repress the subscripts. We denote ${\cal M}(h,f)$ as $hf$. 
       
       The family of morphisms indexed by pairs of objects $(g,f),(g',f') \in 
       \left|{\bcal U}(B,C) \otimes^{(1)}_{1} {\bcal U}(A,B)\right|$ furnishes the direct analogue of horizontal composition of 2-cells
       as can be seen by observing their domain and range in ${\cal V}$:
       $${\cal M}_{ABC_{(g,f)(g',f')}}:[{\bcal U}(B,C) \otimes^{(1)}_{1} 
       {\bcal U}(A,B)]((g,f),(g',f')) \to {\bcal U}(A,C)(gf,g'f')$$
       Recall that 
       $$[{\bcal U}(B,C) \otimes^{(1)}_{1} {\bcal U}(A,B)]((g,f),(g',f')) = {\bcal U}(B,C)(g,g') \otimes_{2} {\bcal U}(A,B)(f,f').$$
       
          \item For each object $A \in \left|{\bcal U}\right|$ a ${\cal V}$--functor
       $${\cal J}_A: {\cal I} \to {\bcal U}(A,A)$$ 
       We denote ${\cal J}_A(0)$ as $1_{A}$. 
          \item (Associativity axiom of a strict ${\cal V}$--2--category.) We require a commuting pentagon.
          Since the morphisms are 
       ${\cal V}$--functors this amounts to saying that the 
       functors given by the two legs of the diagram are equal. 
       For objects we have the equality $(fg)h = f(gh).$  
       
       For the  hom--object morphisms we have the following family of commuting diagrams for associativity, where the first bullet represents
       $$[({\bcal U}(C,D)\otimes^{(1)}_{1} {\bcal U}(B,C)) \otimes^{(1)}_{1} {\bcal U}(A,B)](((f,g),h),((f',g'),h'))$$
       and the reader may fill in the others
       
       $$
        
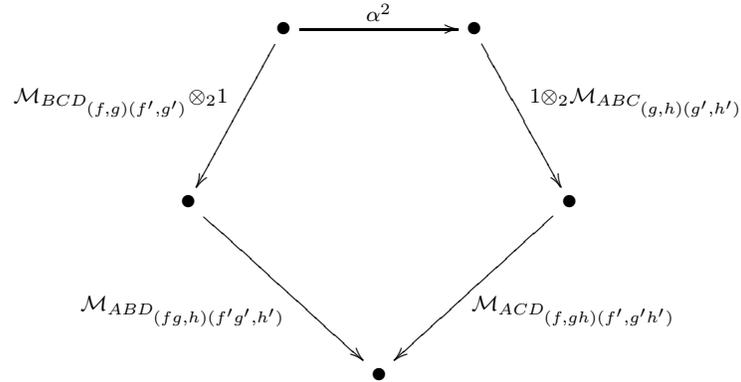
\xymatrix{
         &\bullet
         \ar[rr]^{\alpha^{2}}
         \ar[ddl]_{{\cal M}_{BCD_{(f,g)(f',g')}} \otimes_{2} 1}
         &&\bullet
         \ar[ddr]^{1 \otimes_{2} {\cal M}_{ABC_{(g,h)(g',h')}}}&\\\\
         \bullet
         \ar[ddrr]_{{\cal M}_{ABD_{(fg,h)(f'g',h')}}}
         &&&&\bullet
           \ar[ddll]^{{\cal M}_{ACD_{(f,gh)(f',g'h')}}}
           \\\\&&\bullet&&&
           }$$
      The underlying diagram for this commutativity is
         $$
         
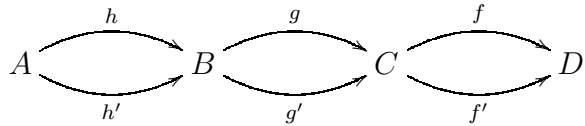
\xymatrix@R-=3pt{
         &&&\\
         A
         \ar@/^1pc/[rr]^h
         \ar@/_1pc/[rr]_{h'}
         &&B
         \ar@/^1pc/[rr]^g
         \ar@/_1pc/[rr]_{g'}
         &&C
         \ar@/^1pc/[rr]^f
         \ar@/_1pc/[rr]_{f'}
         &&D
         \\
         &&&\\
         }
    $$
      \item (Unit axioms of a strict ${\cal V}$--2--category.) We require commuting triangles. 
        For objects we have the equality $f1_{A} = f = 1_{B}f.$
        For the  unit morphisms we have that the triangles in the following diagram commute.
        
        \noindent
		  	      	            	          \begin{center}
		        	          \resizebox{6in}{!}{ 
	  $$
           \xymatrix@C=-5pt{
           [{\cal I}\otimes^{(1)}_1 {\bcal U}(A,B)]((0,f),(0,g))
           \ar[rrd]^{=}
           \ar[dd]_{{\cal J}_{B_{00}}\otimes_2 1}
           &&&&[{\bcal U}(A,B)\otimes^{(1)}_1 {\cal I}]((f,0),(g,0))
           \ar[dd]^{{1}\otimes_2 {\cal J}_{A_{00}}}
           \ar[lld]^{=}\\
           &&{\bcal U}(A,B)(f,g)\\
           [{\bcal U}(B,B)\otimes^{(1)}_1 {\bcal U}(A,B)]((1_B,f),(1_B,g))
           \ar[rru]|{{\cal M}_{ABB_{(1_B,f)(1_B,g)}}}
           &&&&[{\bcal U}(A,B)\otimes^{(1)}_1 {\bcal U}(A,A)]((f,1_A),(g,1_A))
           \ar[llu]|{{\cal M}_{AAB_{(f,1_A)(g,1_A)}}}
           }
        $$ 
        	            }
	  	           	                    \end{center}
     The underlying diagrams for this commutativity are
     
     \noindent
     	  	      	            	          \begin{center}
     	        	          \resizebox{6in}{!}{ 
     	 $$
      \xymatrix@R-=3pt{
      &\ar@{=>}[dd]^{1_{1_A}}&&\\
      A
      \ar@/^1pc/[rr]^{1_A}
      \ar@/_1pc/[rr]_{1_A}
      &&A
      \ar@/^1pc/[rr]^f
      \ar@/_1pc/[rr]_g
      &&B
      \\
      &&&&\text{ }&\\
      }
      =
      \xymatrix@R-=3pt{
           &&\\
          A
          \ar@/^1pc/[rr]^f
          \ar@/_1pc/[rr]_g
          &&B
         \\
        &&&
        }
      =
      \xymatrix@R-=3pt{
        &&&\ar@{=>}[dd]^{1_{1_B}}\\
        A
        \ar@/^1pc/[rr]^f
        \ar@/_1pc/[rr]_g
        &&B
        \ar@/^1pc/[rr]^{1_B}
        \ar@/_1pc/[rr]_{1_B}
        &&B
        \\
        &&&\\
      }
    $$ 
                }
    	  	           	                    \end{center}
    
    \end{enumerate}
         
     \begin{theorem}   Consequences of ${\cal V}$--functoriality of ${\cal M}$ and ${\cal J}$:
       First the ${\cal V}$--functoriality of ${\cal M}$ implies that the following (expanded) diagram commutes
       
       \noindent
       		          \begin{center}
	          \resizebox{6.5in}{!}{
       \begin{footnotesize}
       $$
        \xymatrix@C=-140pt{
        &({\bcal U}(B,C)(k,m)\otimes_1 {\bcal U}(B,C)(h,k))\otimes_2 ({\bcal U}(A,B)(g,l)\otimes_1 {\bcal U}(A,B)(f,g))
       \ar[rdd]^{M_{hkm}\otimes_2 M_{fgl}}
       \\\\
        ({\bcal U}(B,C)(k,m)\otimes_2 {\bcal U}(A,B)(g,l))\otimes_1 ({\bcal U}(B,C)(h,k)\otimes_2 {\bcal U}(A,B)(f,g))
       \ar[dd]^{{\cal M}_{ABC_{(k,g)(m,l)}}\otimes_1 {\cal M}_{ABC_{(h,f)(k,g)}}}
       \ar[ruu]^{\eta^{1,2}}
       &&{\bcal U}(B,C)(h,m)\otimes_2 {\bcal U}(A,B)(f,l)
       \ar[dd]^{{\cal M}_{ABC_{(h,f)(m,l)}}}
       \\\\
       {\bcal U}(A,C)(kg,ml)\otimes_1 {\bcal U}(A,C)(hf,kg)
       \ar[rr]^{M_{(hf)(kg)(ml)}}
       &&{\bcal U}(A,C)(hf,ml)
       \\
       }
       $$
       \end{footnotesize}
       }
       	                    \end{center}
	                    
       The underlying diagram is
    $$
    \xymatrix@R-=16pt{
    &&&
    \\
    A
    \ar@/^2pc/[rr]^f
    \ar[rr]^g
    \ar@/_2pc/[rr]^l
    && B
    \ar@/^2pc/[rr]^h
    \ar[rr]^k
    \ar@/_2pc/[rr]^m
    && C\\
    &&&&&
    }
    $$       
       
       Secondly the ${\cal V}$--functoriality of ${\cal M}$ implies that the following (expanded) diagram commutes
       $$
         \xymatrix{
         &&{\bcal U}(B,C)(g,g)\otimes_2 {\bcal U}(A,B)(f,f)
         \ar[dd]^{{\cal M}_{ABC_{(g,f)(g,f)}}}\\
         I
         \ar[rru]^{j_{g}\otimes_2 j_{f}}
         \ar[rrd]_{j_{gf}}\\
         &&{\bcal U}(A,C)(gf,gf)
         }
        $$
     The underlying diagram here is 
    $$
    \xymatrix@R-=3pt{
    &\ar@{=>}[dd]^{1_f}&&\ar@{=>}[dd]^{1_g}\\
    A
    \ar@/^1pc/[rr]^f
    \ar@/_1pc/[rr]_f
    &&B
    \ar@/^1pc/[rr]^g
    \ar@/_1pc/[rr]_g
    &&C
    \\
    &&&\\
    }
    =
    \xymatrix@R-=3pt{
    &\ar@{=>}[dd]^{1_{gf}}&&\\
    A
    \ar@/^1pc/[rr]^{gf}
    \ar@/_1pc/[rr]_{gf}
    &&C
    \\
    &&&\\
    } 
    $$   
       
       In addition, the ${\cal V}$--functoriality of ${\cal J}$ implies that the following (expanded) diagram commutes
       $$
         \xymatrix{
         &&{\cal I}(0,0)
         \ar[dd]^{{\cal J}_{A_{00}}}\\
         I
         \ar[rru]^{j_{0}}
         \ar[rrd]_{j_{1_A}}\\
         &&{\bcal U}(A,A)(1_A,1_A)
         }
       $$
       Which means that 
       $${\cal J}_{A_{00}}: I \to {\bcal U}(A,A)(1_{A},1_{A}) = j_{1_{A}}.$$
 \end{theorem}       
   
\newpage   
\MySection{Category of ${\cal V}$--$n$--Categories}    
      The definition of a category enriched over ${\cal V}$--$(n-1)$--Cat is simply
        stated by describing the process as enriching over 
        ${\cal V}$--$(n-1)$--Cat with the first of the $k-n$ ordered products. In  detail this means that:
        
        \begin{definition}  \label{super} A (small, strict) ${\cal V}$--{\it n--category}
        ${\bcal U}$ consists of
                
                \begin{enumerate}
                    \item A set of objects $\left|{\bcal U}\right|$
                    \item For each pair of objects $A,B \in \left|{\bcal U}\right|$ a
                    ${\cal V}$--$(n-1)$--category ${\bcal U}(A,B).$
                    \item For each triple of objects $A,B,C \in 
                \left|{\bcal U}\right|$ a ${\cal V}$--$(n-1)$--functor 
                $${\bcal M}_{ABC}:{\bcal U}(B,C) \otimes^{(n-1)}_{1} {\bcal U}(A,B) \to {\bcal U}(A,C)$$
                    \item For each object $A \in \left|{\bcal U}\right|$ a ${\cal V}$--$(n-1)$--functor
                $${\bcal J}_A: {\bcal I}^{(n-1)} \to {\bcal U}(A,A)$$ 
                Henceforth we let the dimensions of domain for and particular instances of ${\bcal M}$ and ${\bcal J}$
                largely be determined by context.
                    \item Axioms: The ${\cal V}$--$(n-1)$--functors that play the
        role of composition and identity
                    obey commutativity of a pentagonal diagram (associativity
        axiom) and of two triangular diagrams (unit axioms).
                   This amounts to saying that the 
                functors given by the two legs of each diagram are equal. 
                $$
      	          \xymatrix@R=10pt@C=10pt{
      	           &\bullet
      	           \ar[rr]^{\alpha^{(n)}}
      	           \ar[ddl]_{{\bcal M}_{BCD} \otimes^{(n)}_{1} 1}
      	           &&\bullet
      	           \ar[ddr]^{1 \otimes^{(n)}_{1} {\bcal M}_{ABC}}&\\\\
      	           \bullet
      	           \ar[ddrr]_{{\bcal M}_{ABD}}
      	           &&&&\bullet
      	             \ar[ddll]^{{\bcal M}_{ACD}}
      	             \\\\&&\bullet&&&
                 }$$
                 $$
      	              \xymatrix@C=-5pt{
      	              {\bcal I}^{(n)}\otimes^{(n)}_1 {\bcal U}(A,B)
      	              \ar[rrd]^{=}
      	              \ar[dd]_{{\bcal J}_{B}\otimes^{(n)}_1 1}
      	              &&&&{\bcal U}(A,B)\otimes^{(n)}_1 {\bcal I}^{(n)}
      	              \ar[dd]^{{1}\otimes^{(n)}_1 {\bcal J}_{A}}
      	              \ar[lld]^{=}\\
      	              &&{\bcal U}(A,B)\\
      	              \bullet
      	              \ar[rru]|{{\bcal M}_{ABB}}
      	              &&&&\bullet
      	              \ar[llu]|{{\bcal M}_{AAB}}
      	              }
              $$
              The consequences of these axioms are expanded commuting diagrams just as in our example for $n=2.$ 
                 \end{enumerate}
       \end{definition}      
        
        This definition requires that there be definitions of the unit ${\bcal I}^{(n)}$ and of ${\cal
        V}$--$n$--functors in place. First, from the proof of monoidal structure on ${\cal
        V}$--$n$--Cat, we can infer a recursively defined 
        unit ${\cal V}$--$n$--category.
        
        \begin{definition} The unit object in ${\cal V}$--$n$--Cat is the ${\cal
        V}$--$n$--category 
          ${\bcal I}^{(n)}$ with one object $\mathbf{0}$ and with ${\bcal
        I}^{(n)}(\mathbf{0},\mathbf{0}) = {\bcal I}^{(n-1)},$
          where ${\bcal I}^{(n-1)}$ is the unit object in ${\cal
        V}$--$(n-1)$--Cat. Of course we let ${\bcal I}^{(0)}$ be $I$ in
          ${\cal V}.$ Also ${\bcal M}_{000} = {\bcal J}_{0} = 1_{{\bcal I^{(n)}}}.$
          \end{definition}
          
        Now we can define the functors:  
        
                  \begin{definition} \label{onecell}
                  For two ${\cal V}$--$n$--categories ${\bcal U}$ and ${\bcal W}$
        a ${\cal V}${\it --$n$--functor} $T:{\bcal U}\to{\bcal W}$ is
        a function on objects 
                  $\left|{\bcal U}\right|\to\left|{\bcal W}\right|$ and a
        family of ${\cal V}$--$(n-1)$--functors
                  $T_{UU'}:{\bcal U}(U,U')\to{\bcal W}(TU,TU').$ These latter
        obey commutativity of the usual diagrams.
                  \begin{enumerate}
                   \item For $U,U',U'' \in \left|{\bcal U}\right|$
                   $$
                     \xymatrix@C=45pt@R=25pt{
                     &\bullet
                     \ar[rr]^{{\bcal M}_{UU'U''}}
                     \ar[d]^{T_{U'U''} \otimes^{(n-1)}_1 T_{UU'}}
                     &&\bullet
                     \ar[d]^{T_{UU''}}&\\
                     &\bullet
                     \ar[rr]_{{\bcal M}_{(TU)(TU')(TU'')}}
                     &&\bullet
                     }
                   $$
                   \item
                   $$
                     \xymatrix{
                     &&\bullet
                     \ar[dd]^{T_{UU}}\\
                     {\bcal I}^{(n-1)}
                     \ar[rru]^{{\bcal J}_U}
                     \ar[rrd]_{{\bcal J}_{TU}}\\
                     &&\bullet
                     }
                   $$
                  Here a ${\cal V}$--$0$--functor is just a morphism in ${\cal V}.$
                 \end{enumerate}
                  \end{definition}
                  
          ${\cal V}$--$n$--categories and ${\cal V}$--$n$--functors form a category. 
          Composition of ${\cal V}$--$n$--functors is just composition of the object functions and 
	    composition of the hom--category ${\cal V}$--$(n-1)$--functors, with appropriate subscripts. Thus
	    $(ST)_{UU'}(f) = S_{TUTU'}(T_{UU'}(f)).$ Then it is straightforward to verify that the axioms are obeyed,
	    as in $$(ST)_{U'U''}(f)(ST)_{UU'}(g)$$ 
	    $$= S_{TU'TU''}(T_{U'U''}(f))S_{TUTU'}(T_{UU'}(g))$$ 
	    $$= S_{TUTU''}(T_{U'U''}(f)T_{UU'}(g))$$
	    $$= S_{TUTU''}(T_{UU''}(fg))$$
	    $$= (ST)_{UU''}(fg).$$
	    That this composition is associative follows from the associativity of composition of the underlying functions and
	    ${\cal V}$--$(n-1)$--functors. The 2--sided identity for this composition $1_{\bcal U}$ is made of the 
	    identity function (on objects) and
  identity ${\cal V}$--$(n-1)$--functor (for hom--categories.)
               The 1--cells we have just defined play a special role in the  
    definition of a general $k$--cell for $k\ge 2$. These higher morphisms will be shown to exist 
and described in some detail in section 5.

\newpage
 \MySection{Induced Functors}         
  
   Recall that in \cite{forcey1} we proved inductively that ${\cal V}$--$n$--Cat is indeed a $(k-n)$--fold monoidal category.
Next we describe a similar result for $k$--fold monoidal functors. This type of fact is often labeled a change of base
theorem. First the definitions of monoidal and $k$--fold monoidal functors, taken from \cite{forcey1} which in turn heavily 
utilized \cite{Balt}. 

\begin{definition}
    A {\it monoidal functor} $(F,\eta) :({\cal C},I)\to({\cal D},J)$ between monoidal categories consists
    of a functor $F$ such that $F(I)=J$ together with a natural transformation
    $$
    \eta_{AB}:F(A)\otimes F(B)\to F(A\otimes B),
    $$
    which satisfies the following conditions
    \begin{enumerate}
    \item Internal Associativity: The following diagram commutes 
    $$
    \diagram
    (F(A)\otimes F(B))\otimes F(C)
    \rrto^{\eta_{AB}\otimes 1_{F(C)}}
    \dto^{\alpha}
    &&F(A\otimes B)\otimes F(C)
    \dto^{\eta_{(A\otimes B)C}}\\
    F(A)\otimes (F(B)\otimes F(C))
    \ar[d]^{1_{F(A)}\otimes \eta_{BC}}
    &&F((A\otimes B)\otimes C)
    \ar[d]^{F\alpha}\\
    F(A)\otimes F(B\otimes C)
    \rrto^{\eta_{A(B\otimes C)}}
    &&F(A\otimes (B\otimes C))
    \enddiagram
    $$
    \item Internal Unit Conditions: $\eta_{AI}=\eta_{IA}=1_{F(A)}.$
    \end{enumerate}
    \end{definition}
    Given two monoidal functors $(F,\eta) :{\cal C}\to{\cal D}$ and $(G,\zeta) 
    :{\cal D}\to{\cal E}$,
    we define their composite to be the monoidal functor $(GF,\xi) :
    {\cal C}\to{\cal E}$, where
    $\xi$ denotes the composite
    $$
    \diagram
    GF(A)\otimes GF(B)\rrto^{\zeta_{F(A)F(B)}} 
    && G\bigl(F(A)\otimes F(B)\bigr)\rrto^{G(\eta_{AB})}
    &&GF(A\otimes B).
    \enddiagram
    $$
    It is easy to verify that $\xi$ satisfies the internal associativity condition above by subdividing the 
    necessary commuting diagram into two regions that commute by the axioms for $\eta$ and $\zeta$ respectively 
    and two that commute due to their naturality.

\begin{definition} An {\it $n$--fold monoidal functor}
  $(F,\lambda^1,\dots,\lambda^n):{\cal C}\to{\cal D}$ between $n$--fold monoidal categories
  consists of a functor $F$ such that $F(I)=J$ together with natural
  transformations
  $$\lambda^i_{AB}:F(A)\otimes_i F(B)\to F(A\otimes_i B)\quad i=1,2,\dots, n$$
  satisfying the same associativity and unit conditions as monoidal functors.
  In addition the following hexagonal interchange diagram commutes:
  $$
  \diagram
  (F(A)\otimes_j F(B))\otimes_i(F(C)\otimes_j F(D))
  \xto[rrr]^{ {\eta^{ij}_{F(A)F(B)F(C)F(D)}}}
  \dto^{ {\lambda^j_{AB}\otimes_i\lambda^j_{CD}}}
  &&&(F(A)\otimes_i F(C))\otimes_j(F(B)\otimes_i F(D))
  \dto^{ \lambda^i_{AC}\otimes_j\lambda^i_{BD}}\\
  F(A\otimes_j B)\otimes_i F(C\otimes_j D)
  \dto^{ {\lambda^i_{(A\otimes_j B)(C\otimes_j D)}}}
  &&&F(A\otimes_i C)\otimes_j F(B\otimes_i D)
  \dto^{ {\lambda^j_{(A\otimes_i C)(B\otimes_i D)}}}\\
  F((A\otimes_j B)\otimes_i(C\otimes_j D))
  \xto[rrr]^{ F(\eta^{ij}_{ABCD})}
  &&& F((A\otimes_i C)\otimes_j(B\otimes_i D))
  \enddiagram
  $$
  \end{definition}
  Composition of $n$-fold monoidal functors is defined 
  as for monoidal functors.  

The following generalizes a well known theorem for the symmetric case, found for instance in \cite{Borc}.
\begin{theorem}\label{cob} Change of Base:
Suppose we have an $n$--fold monoidal functor  $(F,\lambda^i):{\cal V}\to{\cal W}$. 
Then there exists an $(n-1)$--fold monoidal functor $(F^{(1)},\lambda^{(1)i}):{\cal V}$--Cat$\to {\cal W}$--Cat.
\end{theorem}

{\bf Proof } 
Let ${\cal A}$ and ${\cal B}$ be ${\cal V}$--categories and $G:{\cal A}\to{\cal B}$ a ${\cal V}$--functor between them.
First we define the action of $F^{(1)}$ on objects by describing the ${\cal W}$--category $F^{(1)}({\cal A}).$ We define
$\left|F^{(1)}({\cal A})\right| = \left|{\cal A}\right|$ and
for $A, A'$ objects in ${\cal A}$ we let $F^{(1)}({\cal A})(A,A')=F({\cal A}(A,A'))$. For a triple of objects $A,A',A''$ in 
$\left|{\cal A}\right|$ we have the composition morphism $M$ in ${\cal W}$ for $F^{(1)}({\cal A})$ given by 

$$
\xymatrix{
F^{(1)}({\cal A})(A',A'')\otimes_1 F^{(1)}({\cal A})(A,A')
\ar@{=}[d]
\\
F({\cal A}(A',A''))\otimes_1 F({\cal A}(A,A'))
\ar[d]^{\lambda^1}
\\
F({\cal A}(A',A'')\otimes_1 {\cal A}(A,A'))
\ar[d]^{F(M_{AA'A''})}
\\
F({\cal A}(A,A''))
\ar@{=}[d]
\\
F^{(1)}({\cal A})(A,A')
}
$$
 and the unit morphism $j_A:J\to F^{(1)}({\cal A})(A,A)$ is just $F(j_A).$
The ${\cal W}$--functor $F^{(1)}(G)$ has the same underlying function on objects as does $G$ and we define
 $F^{(1)}(G)_{AA'}:F^{(1)}({\cal A})(A,A')\to F^{(1)}({\cal B})(GA,GA')$ to be simply $F(G_{AA'}).$ In \cite{Borc} it is pointed out that 
$F^{(1)}$ is a 2--functor.
The reader should check that these images of $F^{(1)}$ obey the axioms of enriched categories and functors. 

Next we define
$$\lambda^{(1)i}_{{\cal A}{\cal B}}:F^{(1)}({\cal A})\otimes^{(1)}_i F^{(1)}({\cal B})\to F^{(1)}({\cal A}\otimes^{(1)}_i {\cal B})\quad i=1,2,\dots, (n-1)$$
to be the identity on objects and to have hom--object morphisms given by
$$\lambda^{(1)i}_{{{\cal A}{\cal B}}_{(AB)(A'B')}} = \lambda^{(i+1)}_{{\cal A}(A,A'){\cal B}(B,B')}$$
It is clear that the $\lambda^{(1)i}$ will inherit the required internal associativity and unit conditions from 
$\lambda^{(i+1)}.$
It is a useful exercise for the reader to work out how the hexagonal interchange diagram and the naturality of $\lambda^{(i+1)}$
imply the enriched functoriality of $\lambda^{(1)i}.$
 
When attached to a functor, the superscript in parentheses will indicate how many times this inducing has been applied.
In other words, $F^{(n)} = F^{(n-1)(1)}=F^{(1)(n-1)}.$
 We now describe the special case of change of base to that of the category of sets; induced forgetful functors.
\begin{theorem}\label{hom} Given ${\cal V}$ an $n$--fold monoidal category, the functor $\text{Hom}_{\cal V}(I,\_):{\cal V}\to $ {\bf Set} 
is $n$--fold monoidal.
\end {theorem}
 {\bf Proof } 
  The $n$--fold monoidal structure on the (strictly monoidal, symmetric) {\bf Set} is given by $\otimes_1 = \otimes_2 = \dots =\otimes_n= \times.$
  The natural transformations $\lambda^i_{AB}:\text{Hom}_{\cal V}(I,A)\times\text{Hom}_{\cal V}(I,B)\to\text{Hom}_{\cal V}(I,A\otimes_i B)$ are given by
 $\lambda^i_{AB}((f,g))=f\otimes_i g.$ Thus the internal associativity and unit conditions hold trivially.

The action of $\text{Hom}_{\cal V}(I,\_)$ on morphisms is given by $\text{Hom}_{\cal V}(I,f)(h)=f\circ h$ for $h\in \text{Hom}_{\cal V}(I,A)$ and
 $f:A\to B.$
That last fact in mind, it is easy to see that the 
 hexagonal interchange diagrams for $\lambda^i$ commute due to the naturality of $\eta$.

Now a bit of theory that we have been ignoring for simplicity is forced upon us. The unit of {\bf Set} is the single element
set we denote as $\{0\}$, but this is not a strict unit; rather there are canonical isomorphisms $X\times\{0\}\to X; 
\{0\}\times X\to X$ that are natural and obey commuting diagrams as shown in \cite{Borc}. Thus the identification of
$\text{Hom}_{{\cal V}}(I,I)$ with $\{0\}$ is not precise--$\{0\}$ is mapped to $1_I$ and an additional pair of diagrams is shown 
to commute. Again the reader is referred to \cite{Borc} where this part of the proof is given.

The next fact is that when we apply  Theorem~\ref{cob} to the functor of Theorem~\ref{hom} we get not only the 
corollary that there is an induced $(n-1)$--fold monoidal functor 
$\text{Hom}^{(1)}_{{\cal V}}(I,\_):{\cal V}\text{--Cat} \to \text{{\bf Set}--Cat} = \text{ }${\bf Cat},
 but also a direct description of that functor in terms of the unit in the category of enriched categories.

\begin{theorem}\label{hom2}
The following two functors are equivalent:
$$\text{Hom}^{(1)}_{{\cal V}}(I,\_) = \text{Hom}_{({\cal V}\text{--Cat})}({\cal I},\_)$$
where ${\cal I}$ is the unit enriched category in ${\cal V}$--Cat. 

\end{theorem}

{\bf Proof } 
(\cite{Kelly})
The image of ${\cal A}$ under the 
induced functor is known as the underlying category of ${\cal A}$, denoted ${\cal A}_0.$
 It has objects the same as ${\cal A}$ and morphisms ${\cal A}_0(A,A')=\text{Hom}_{\cal V}(I,{\cal A}(A,B)).$ Composition
is given by:
$$
\xymatrix{
I\otimes I
\ar[d]^{f\otimes g}
\\
{\cal A}(A',A'')\otimes{\cal A}(A,A')
\ar[d]^{M}
\\
{\cal A}(A,A'')
}
$$
and the identity in ${\cal A}_0(A,A)$ is $j_A.$ In \cite{Kelly} it is shown that this functor is representable, in fact that 
the equality of the theorem holds. 

Notice that this fact is recursively true for all the unit categories ${\bcal I}^n.$
In other words, we can write 
$\text{Hom}_{{\cal V}\text{--}n\text{--Cat}}^{(k)}({\bcal I}^{(n)},\_)=\text{Hom}_{{\cal V}\text{--}(n+k)\text{--Cat}}({\bcal I}^{(n+k)},\_).$

The point of view taken by the next section differs subtly from that of \cite{Kelly}. In \cite{Kelly} the definition of the 
$2$--cells of ${\cal V}$--Cat is given first, and then the representable 2--functor $(\text{ })_0 = \text{Hom}_{{\cal V}\text{--Cat}}({\cal I},\_)$ is 
described. The author then notes that the image of the enriched natural transformation under the functor in question is an
ordinary natural transformation with the same components as the enriched one. We would like to back up a bit and define the
enriched natural transformations in terms of the image of $(\text{ })_0.$ In other words, let the ${\cal V}$--natural transformations
between ${\cal V}$--functors $T$ and $S$ be the ordinary natural transformations between $T_0$ and $S_0$, with the additional 
requirement that their components obey the commuting diagram of Definition~\ref{enr:nat:trans}. The latter requirement
is strengthening; enriched natural transformations $T\to S$ are automatically ordinary ones $T_0 \to S_0$ but the converse is not true.
This point of view was suggested by the referee as leading to the proper method of showing the general categorical dimension of ${\cal V}$--$n$--Cat.  

\newpage    
\MySection{Categorical Structure of ${\cal V}$--$n$--Cat}
 
 \begin{theorem}\label{general}
 For ${\cal V}$ $k$--fold monoidal the category of ${\cal V}$--$n$--categories has the additional structure of 
 an $(n+1)$--category .
 \end{theorem}
 
 {\bf Proof } 
 This fact results from the theorems of the last section applied to the functor 
 $(\text{ })_{\mathbf 0} = \text{Hom}_{{\cal V}\text{--}n\text{--Cat}}({\bcal I}^{n},\_).$
  We have from the latter's identification with the 
 functor induced by $\text{Hom}_{{\cal V}\text{--}(n-1)\text{--Cat}}({\bcal I}^{(n-1)},\_)$ by induction that its range is
 $n$--Cat. This is due to Theorem~\ref{cob} since $n$--Cat $= (n-1)$--Cat--Cat. 
 
 Now the overall $(n+1)$--category structure of ${\cal V}$--$n$--Cat is derived from the well known $(n+1)$--category 
 structure of $n$--Cat. More precisely, we already have seen that ${\cal V}$--$n$--categories and ${\cal V}$--$n$--functors
 form a category. Let the 2--cells between two ${\cal V}$--$n$--functors $T$ and $S$ be the natural transformations
 between $T_{\mathbf 0}$ and $S_{\mathbf 0}$ and the enriched modifications and higher $k$--cells the ordinary ones that arise on the same principle.
 The axioms of an $(n+1)$--category are thus automatically satisfied by these inherited cells.

Now we demonstrate that 
${\cal V}$--$n$--Cat has as well a special sub--$(n+1)$--category structure that is a restriction of the morphisms occurring 
 between images of the induced forgetful functors. This structure is the natural extension of the 2--category structure 
 described for ${\cal V}$--Cat in \cite{Kelly}, and therefore when we speak of ${\cal V}$--$n$--Cat hereafter it will be the following 
 structure to which we refer.

              \begin{definition} \label{kcell}
	                  A ${\cal V}${\it --n:k--cell} $\alpha$ between $(k-1)$--cells
	        $\psi^{k-1}$ and $\phi^{k-1}$, written 
	                   $$\alpha:\psi^{k-1}\to \phi^{k-1}:\psi^{k-2}\to \phi^{k-2}: ... :\psi^{2}\to \phi^{2}:F\to G:{\bcal U}\to {\bcal W}$$ 
	                  where $F$ and $G$ are  ${\cal V}$--$n$--functors and where the superscripts denote cell dimension,
	                  is a function sending
	                  each $U \in  \left|{\bcal U}\right|$ to a 
	                  ${\cal V}$--$((n-k)+1)$--functor 
	                 $$\alpha_{U}: {\bcal I}^{((n-k)+1)}\to {\bcal W}(FU,GU)(\psi^{2}_U\mathbf{0},\phi^{2}_U\mathbf{0}) ... (\psi^{k-1}_U\mathbf{0},\phi^{k-1}_U\mathbf{0})$$
	                  in such a way that we have commutativity of the following
	        diagram. Note that the final (curved) equal sign is implied recursively by the diagram for the $(k-1)$--cells. 
                  \end{definition}
                  
                  \noindent
	          \begin{center}
	          \resizebox{6.5in}{!}{
	          \begin{footnotesize}
	                     $$
	                      \xymatrix@C=-140pt@R=27pt{
	                      &
	                      &*\txt{${\bcal W}(FU',GU')(\psi^{2}_{U'}\mathbf{0},\phi^{2}_{U'}\mathbf{0}) ... (\psi^{k-1}_{U'}\mathbf{0},\phi^{k-1}_{U'}\mathbf{0})$\\$\otimes^{((n-k)+1)}_{k-1} {\bcal W}(FU,FU')(F(x_2),F(y_2))...(F(x_{k-1}),F(y_{k-1}))$}
	                      \ar[rd]^-{{\bcal M}}
	                    \\
	                    &
	                    {\bcal I}^{((n-k)+1)} \otimes^{((n-k)+1)}_{k-1} {\bcal U}(U,U')(x_2,y_2)...(x_{k-1},y_{k-1})\text{ }\text{ }\text{ }\text{ }\text{ }\text{ }\text{ }\text{ }\text{ }
	                      \ar[ru]^<<{\alpha_{U'} \otimes^{((n-k)+1)}_{k-1} F\text{ }}
	                    &&{\bcal W}(FU,GU')(\psi^{2}_{U'}\mathbf{0}F(x_2),\phi^{2}_{U'}\mathbf{0}F(y_2)) ... (\psi^{k-1}_{U'}\mathbf{0}F(x_{k-1}),\phi^{k-1}_{U'}\mathbf{0}F(y_{k-1})) 
	                    \ar@{=}@/^2pc/[dd]
	                    \\
	                    {\bcal U}(U,U')(x_2,y_2)...(x_{k-1},y_{k-1})
	                      \ar[ru]^{=}
	                    \ar[rd]_{=}
	                    \\
	                    &
	                    {\bcal U}(U,U')(x_2,y_2)...(x_{k-1},y_{k-1}) \otimes^{((n-k)+1)}_{k-1} {\bcal I}^{((n-k)+1)}\text{ }\text{ }\text{ }\text{ }\text{ }\text{ }\text{ }\text{ }\text{ }
	                      \ar[rd]_<<{G \otimes^{((n-k)+1)}_{k-1} \alpha_{U}\text{ }}
	                    &&{\bcal W}(FU,GU')(G(x_2)\psi^{2}_U\mathbf{0},G(y_2)\phi^{2}_U\mathbf{0})...(G(x_{k-1})\psi^{k-1}_U\mathbf{0},G(y_{k-1})\phi^{k-1}_U\mathbf{0})   
	                    \\
	                      &
	                      &*\txt{${\bcal W}(GU,GU')(G(x_2),G(y_2))...(G(x_{k-1}),G(y_{k-1}))$\\$\otimes^{((n-k)+1)}_{k-1} {\bcal W}(FU,GU)(\psi^{2}_U\mathbf{0},\phi^{2}_U\mathbf{0}) ... (\psi^{k-1}_U\mathbf{0},\phi^{k-1}_U\mathbf{0})$}
	                      \ar[ru]_-{{\bcal M}}
	                      }
	                    $$
	                \end{footnotesize}    
	                    }
	                    \end{center}

           Thus for a given value of $n$ there are $k$--cells up to $k = n+1$,  making   ${\cal V}$--$n$--Cat a potential $(n+1)$--category.
           We have already described composition of ${\cal V}$--$n$--functors. Now we describe all other compositions.  
\begin{definition}\label{kcomp} Case 1. 

Let $k=2\dots n+1$ and $i=1\dots k-1.$ Given $\alpha$ and $\beta$ two ${\cal V}$--{\it n:k}--cells that share a common 
${\cal V}$--{\it n:(k-i)}--cell $\gamma,$
  we can compose along the latter morphism as follows
          
          \noindent
	  	          \begin{center}
	          \resizebox{6.5in}{!}{
           $$
           \xymatrix@R=3pt{
           (\beta \circ \alpha)_U =
           {\bcal I}^{((n-k)+1)} = {\bcal I}^{((n-k)+1)} \otimes^{((n-k)+1)}_{i} {\bcal I}^{((n-k)+1)}\text{ }\text{ }\text{ }\text{ }\text{ }\text{ }\text{ }       
           \ar[ddd]^{\beta_U \otimes^{((n-k)+1)}_{i} \alpha_U}
           \\\\\\
           \text{\resizebox{4in}{!}{${\bcal W}(FU,GU)(\psi^{2}_U\mathbf{0},\phi^{2}_U\mathbf{0}) ... (\psi^{k-i-1}_U\mathbf{0},\phi^{k-i-1}_U\mathbf{0})(\gamma_U\mathbf{0},\gamma''_U\mathbf{0})(\psi^{k-i+1}_U\mathbf{0},\phi^{k-i+1}_U\mathbf{0}) ... (\psi^{k-1}_U\mathbf{0},\phi^{k-1}_U\mathbf{0})$}}
           \\
           *\txt{\resizebox{4.25in}{!}{${\scriptscriptstyle\otimes^{((n-k)+1)}_{i} {\bcal W}(FU,GU)(\psi^{2}_U\mathbf{0},\phi^{2}_U\mathbf{0}) ... (\psi^{k-i-1}_U\mathbf{0},\phi^{k-i-1}_U\mathbf{0})(\gamma'_U\mathbf{0},\gamma_U\mathbf{0})(\delta^{k-i+1}_U\mathbf{0},\xi^{k-i+1}_U\mathbf{0}) ... (\delta^{k-1}_U\mathbf{0},\xi^{k-1}_U\mathbf{0})}$}}
           \ar[ddd]^{{\bcal M}}
           \\\\\\
           *\txt{\resizebox{4.5in}{!}{${\scriptscriptstyle{\bcal W}(FU,GU)(\psi^{2}_U\mathbf{0},\phi^{2}_U\mathbf{0}) ... (\psi^{k-i-1}_U\mathbf{0},\phi^{k-i-1}_U\mathbf{0})(\gamma'_U\mathbf{0},\gamma''_U\mathbf{0})(\psi^{k-i+1}_U\mathbf{0}\delta^{k-i+1}_U\mathbf{0},\phi^{k-i+1}_U\mathbf{0}\xi^{k-i+1}_U\mathbf{0}) ... (\psi^{k-1}_U\mathbf{0}\delta^{k-1}_U\mathbf{0},\phi^{k-1}_U\mathbf{0}\xi^{k-1}_U\mathbf{0})}$}}
           }
           $$
           }
	                    \end{center}
	                    \end{definition}
For $\alpha$ and $\beta$ of different dimension and sharing a common cell of dimension lower than either the composition
is accomplished by first raising the dimension of the lower of $\alpha$ and $\beta$ to match the other by replacing it with
a unit (see next Definition.)

{\it
Case 2.

It remains to describe composing along a 0--cell, i.e. along a common ${\cal V}$--$n$--category ${\bcal W}$. We describe composing a higher enriched
cell with an enriched functor, and then leave the remaining possibilities to be accomplished by applying the first case
to the results of such whiskering. 

Composition with a ${\cal V}$--$n$--functor $K:{\bcal W}\to{\bcal X}$ on the right is given by:
$$\xymatrix{(K\alpha)_{U}} = \xymatrix{
                           {\bcal I}^{(n-k)+1)}  
                           \ar[d]_{\alpha_{U}}
                           \\{\bcal W}(FU,GU)(\psi^{2}_U\mathbf{0},\phi^{2}_U\mathbf{0}) ... (\psi^{k-1}_U\mathbf{0},\phi^{k-1}_U\mathbf{0})
                            \ar[d]_{K_{FU,GU}}
                           \\{\bcal X}(KFU,KGU)(K\psi^{2}_U\mathbf{0},K\phi^{2}_U\mathbf{0}) ... (K\psi^{k-1}_U\mathbf{0},K\phi^{k-1}_U\mathbf{0})
                            }$$
                                  
Composing with a ${\cal V}$--$n$--functor $H:{\bcal V}\to{\bcal U}$ on the left is given by $(\alpha H)_V = \alpha_{HV}.$  
}


We describe unit $k$--cells for the above compositions.
\begin{definition}
A unit ${\cal V}${\it --n:k--cell} ${\mathbf 1}_{\psi^{k-1}}$ from  a $(k-1)$--cell
	        $\psi^{k-1}$ to itself,  sends
	                  each $U \in  \left|{\bcal U}\right|$ to 
	                  the ${\cal V}$--$((n-k)+1)$--functor 
	                 $${\bcal J}_{\psi^{k-1}_U\mathbf{0}}: {\bcal I}^{((n-k)+1)}\to {\bcal W}(FU,GU)(\psi^{2}_U\mathbf{0},\phi^{2}_U\mathbf{0}) ... (\psi^{k-1}_U\mathbf{0},\psi^{k-1}_U\mathbf{0})$$
	                  It is straightforward to verify that these fulfill the requirements of Definition~\ref{kcell} and 
	                  indeed are units with respect to Definition~\ref{kcomp}. 
	                  
\end{definition}
Of course the unit for composition along 
	                  a common cell of dimension more than 1 less than the composed cells is constructed of units for all the dimensions
	                  between that of the composed cells and that of the common cell. For example, the unit for composing along a common
	                  0--cell may appear as follows:
	                  $$
			    \xymatrix@R=10pt@C=20pt{
			                   &&\ar@/_.5pc/@{=>}[dd]_{1_{1_{\bcal U}}}
			                   &&\ar@/^.5pc/@{=>}[dd]^{1_{1_{\bcal U}}}
			                   \\
			                   {\bcal U}
			                   \ar@/^3pc/[rrrrrr]^{1_{\bcal U}}
			                   \ar@/_3pc/[rrrrrr]_{1_{\bcal U}}
			                   &&
			                   \ar@3{->}[rr]_{1_{1_{1_{\bcal U}}}}
			                   &&&&{\bcal U}
			                   \\
			                   &&&&&&&&&&
			                   }
  $$    
	                  
This last string of definitions is best grasped by looking at examples. 
The cases for $n = 1$  and for $n = 2$ are
treated carefully in \cite{forcey2}, where the following theorem is proven by brute force for $n=2.$

 \begin{theorem}\label{specific}
 ${\cal V}$--$n$--categories, ${\cal V}$--$n$--functors, and ${\cal V}$--\it{n:k}--cells for $k=2\dots n+1$ together 
have the structure of an $(n+1)$--category.
 \end{theorem} 
  
 {\bf Proof } 
 Two basic sets of facts need checking before we can take advantage of the structure shown in Theorem~\ref{general} First, we need to check that the ${\cal V}$--{\it n:k}--cells and 
 their compositions defined above are indeed a subset of
 the ones which we showed exist between underlying functors in the previous theorem, and ordinary compositions in $n$--Cat
 respectively. Second, we need to check that the ${\cal V}$--{\it n:k}--cells are closed under those compositions. 
Both of these truths are in part shown by use of an induction on $n$.
Note that since the index $n$ is found only in product superscripts and unit categories in Definition~\ref{kcell}
then if we assume that ${\cal V}$--$(n-1)$--Cat (with the $k$--cells of Definition~\ref{kcell}) is an $n$--category
we have the implication that  ${\cal V}$--$n$--Cat is also an $n$--category. The $k$--cells for $k=1,\dots, n$ of 
${\cal V}$--$n$--Cat are defined in exactly the same way as for ${\cal V}$--$(n-1)$--Cat, just with the index $n$
increased. So to complete the inductive step, all that needs to be dealt
with are the top--dimensional cells; the $(n+1)$--cells of ${\cal V}$--$n$--Cat. 
 
 We need to check then that the said ${\cal V}$--{\it n:}$(n+1)$--cells are indeed $(n+1)$--cells in $n$--Cat between the underlying
 $n$--categories ${\bcal U}_{\mathbf 0}$ and ${\bcal W}_{\mathbf 0}$ of ${\bcal U}$ and ${\bcal W}.$
 $(n+1)$--cells in $n$--Cat are collections: for each $U\in \left|{\bcal U}_{\mathbf 0}\right|$ an $n$--cell
 in ${\bcal W}_{\mathbf 0}.$ 
 
 Recall that 
 from Theorem~\ref{hom2}  ${\bcal W}_{\mathbf 0} = \text{Hom}_{{\cal V}\text{--}(n-1)\text{--Cat}}^{(1)}({\bcal I}^{(n-1)},\_)({\bcal W})
 = \text{Hom}_{{\cal V}\text{--}n\text{--Cat}}({\bcal I}^{(n)},{\bcal W})$ 
 
 Our inductive assumption is that ${\cal V}$--$(n-1)$--Cat is an $n$--category. 
 Thus ${\bcal W}_{\mathbf 0}(A,B) = 
 \text{Hom}_{{\cal V}\text{--}(n-1)\text{--Cat}}({\bcal I}^{(n-1)},{\bcal W}(A,B))$ is an $(n-1)$--category itself. 
 Repeated use of Theorem~\ref{hom2} on this hom--$(n-1)$--category yields a complete
 characterization of ${\bcal W}_{\mathbf 0}$.
  From this point of view objects of ${\bcal W}_{\mathbf 0}$ are those of ${\bcal W}$,  1--cells
are objects in $\left|{\bcal W}(A,B)\right|$ for $A,B\in \left|{\bcal W}\right|$, 2--cells are objects in
$\left|{\bcal W}(A,B)(f,g)\right|$ and so on until the highest order cells: $n$--cells are morphisms in ${\cal V}$,
elements of $\text{Hom}_{{\cal V}}(I,{\bcal W}(A,B)(f_2,g_2)\dots(f_n,g_n)).$  So for $k=n+1$ the $k$--cells of Definition~\ref{kcell} 
each are made up of just the right sort of elements to qualify to be $n+1$--cells in $n$--Cat.
Furthermore, the composition of these $n$--cells in ${\bcal W}_{\mathbf 0}$, which is the basis of component--wise composition of the ordinary
$(n+1)$--cells in $n$-Cat, is just as given in Definition~\ref{kcomp}. We can say so simply because that is how the composition in the underlying
categories comes to be defined in the proof of Theorem~\ref{cob}. This in mind, we note that the condition of the commuting diagram
in Definition~\ref{kcell} implies the higher naturality condition satisfied by an $n$--cell in $n$--Cat. 

The second part of our induction concerns the closure of
 the enriched $(n+1)$--cells under their compositions with each other and
with smaller dimensional cells. Given that we have the structure of an $(n+1)$--category
we can assume that the latter sort of composition, generally known as whiskering, is described in terms of the former
by composing with appropriate unit $(n+1)$--cells. Thus we need only check for closure the various cases involved in 
composing two ${\cal V}$--{\it n:}$(n+1)$--cells along common cells of lower dimension $m=k-i=n+1-i$. There are three cases:
$n\ge m\ge 2$, $m=1$, and $m=0.$
For the first case we have the following diagram which shows that the result of  the composition $\nu \circ \mu$ given by
Definition~\ref{kcomp} obeys the axiom given in Definition~\ref{kcell}.

The first bullet in the following diagram is ${\bcal U}(U,U')(x_2,y_2)...(x_{n},y_{n}).$
Let
\noindent
	  	          \begin{center}
\resizebox{6.5in}{!}{
$$\xymatrix@R=3pt{\\X=}\xymatrix@R=-1pt{\text{{\tiny $n$ total sets of parentheses}}\\\text{\resizebox{3in}{!}{${\bcal U}\overbrace{(U,U')(x_2,y_2)...(x_{m-1},y_{m-1})(x_{m},x_{m})(1_{x_{m}},1_{x_{m}})...\left(1_{1_{\ddots_{1_{x_{m}}}}},1_{1_{\ddots_{1_{x_{m}}}}}\right)}$}}}$$
}
\end{center}
  Other vertices in the diagram include:
  $$A= (I\otimes_{i} I) \otimes_n {\bcal U}(U,U')(x_2,y_2)...(x_{n},y_{n})\text{ }$$
  $$B= (I\otimes_{i} I) \otimes_n ({\bcal U}(U,U')(x_2,y_2)...(x_{n},y_{n})\otimes_{i} I)$$
  $$C= (I\otimes_{i} I) \otimes_n ({\bcal U}(U,U')(x_2,y_2)...(x_{n},y_{n})\otimes_{i} X$$ 
  $$D= (I\otimes_n {\bcal U}(U,U')(x_2,y_2)...(x_{n},y_{n})) \otimes_{i} (I\otimes_n X)$$
  $$E= {\bcal U}(U,U')(x_2,y_2)...(x_{n},y_{n})\otimes_{i} I$$
  $$H= {\bcal U}(U,U')(x_2,y_2)...(x_{n},y_{n})\otimes_{i} X$$
  $$K= {\bcal U}(U,U')(x_2,y_2)...(x_{n},y_{n})\otimes_n (I\otimes_{i} I)\text{ }$$  
  $$L= ({\bcal U}(U,U')(x_2,y_2)...(x_{n},y_{n})\otimes_{i} I)\otimes_n (I\otimes_{i} I)$$
  $$N= ({\bcal U}(U,U')(x_2,y_2)...(x_{n},y_{n})\otimes_n I) \otimes_{i} (X\otimes_n I)\text{ }$$
  $$P= ({\bcal U}(U,U')(x_2,y_2)...(x_{n},y_{n})\otimes_n X) \otimes_{i} (I\otimes_n I).$$
  \clearpage
   \noindent
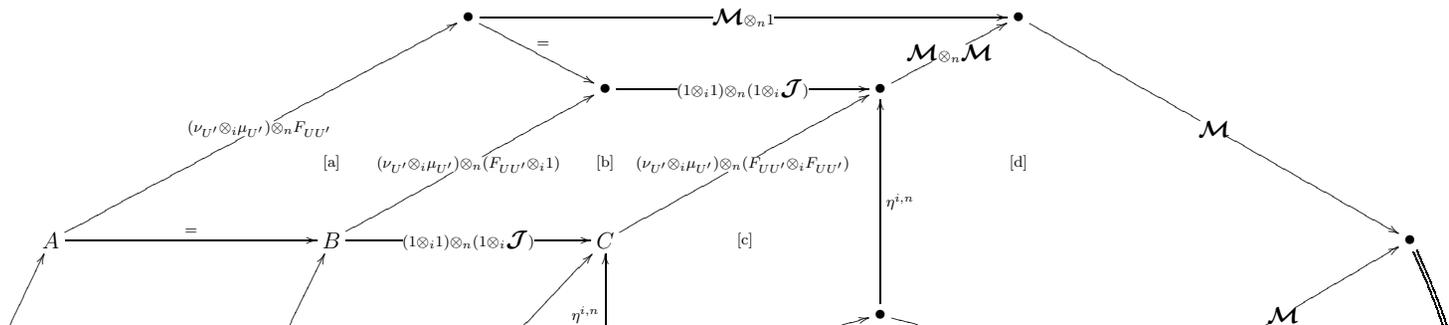
   \begin{center}
	          \resizebox{5.5in}{!}{
   \begin{sideways}
    $$
    \xymatrix@C=57pt@R=26pt{
    &&&&\bullet
    \ar[dr]^{=}
    \ar[rrrr]|{{\bcal M}\otimes_n 1}
    &&&&\bullet
    \ar[dddrrr]|{{\bcal M}}
    \\
    &&&&&\bullet
    \ar[rr]|{(1\otimes_{i} 1)\otimes_n (1\otimes_{i} {\bcal J})}
    &&\bullet
    \ar[ur]|{{\bcal M}\otimes_n {\bcal M}}
    \\
    &&&\text{\scriptsize {[a]}}
    &&\text{\scriptsize {[b]}}
    &&&\text{\scriptsize {[d]}}
    \\
    &A
    \ar[uuurrr]|{(\nu_{U'}\otimes_{i} \mu_{U'})\otimes_n F_{UU'}}
    \ar[rr]^{=}
    &&B
    \ar[uurr]|{(\nu_{U'}\otimes_{i} \mu_{U'})\otimes_n (F_{UU'}\otimes_{i} 1)}
    \ar[rr]|{(1\otimes_{i} 1)\otimes_n (1\otimes_{i} {\bcal J})}
    &&C
    \ar[uurr]|{(\nu_{U'}\otimes_{i} \mu_{U'})\otimes_n (F_{UU'}\otimes_{i} F_{UU'})}
    &\text{\scriptsize {[c]}}
    &
    &&&&\bullet
    \ar@{=}@/^3pc/[dddddddd]
    \\
    &&&&&&&\bullet
    \ar[uuu]_{\eta^{i,n}}
    \ar[drr]|{{\bcal M}\otimes_{i} {\bcal M}}
    \\
    &\text{\scriptsize {[e]}}
    &&\text{\scriptsize {[f]}}
    &&D
    \ar[uu]^{\eta^{i,n}}
    \ar[urr]|{(\nu_{U'}\otimes_n F_{UU'})\otimes_{i} ( \mu_{U'}\otimes_n F_{UU'})}
    &&&&\bullet
    \ar[uurr]|{{\bcal M}}
    \ar@{=}@/^2pc/[dddd]
    \\\\
    \bullet
    \ar[uuuur]^{=}
    \ar[rr]^{=}
    \ar[ddddr]^{=}
    &&E
    \ar[uuuur]^{=}
    \ar[r]|{1\otimes_{i} {\bcal J}}
    \ar[ddddr]^{=}
    &H
    \ar[uurr]^{=}
    \ar[ddrr]^{=}
    \ar[uuuurr]^{=}
    \ar[ddddrr]^{=}
    &&&\text{\scriptsize {[g]}}
    \\\\
    &\text{\scriptsize {[h]}}
    &&\text{\scriptsize {[i]}}
    &&N
    \ar[dd]_{\eta^{i,n}}
    \ar[drr]|{(G_{UU'}\otimes_n \nu_{U})\otimes_{i} (G_{UU'}\otimes_n \mu_{U})}
    &&&&\bullet
    \ar[ddrr]|{{\bcal M}}
    \\
    &&&&&&&\bullet
      \ar[ddd]^{\eta^{i,n}}
    \ar[urr]|{{\bcal M}\otimes_{i} {\bcal M}}
    \\
    &K
    \ar[dddrrr]|{G_{UU'}\otimes_n (\nu_U\otimes_{i} \mu_U)}
    \ar[rr]^{=}
    &&L
    \ar[ddrr]|{(G_{UU'}\otimes_{i} 1)\otimes_n (\nu_U\otimes_{i} \mu_U)}
    \ar[rr]|{(1\otimes_{i} {\bcal J})\otimes_n (1\otimes_{i} 1)}
    &&P
    \ar[ddrr]|{(G_{UU'}\otimes_{i} G_{UU'})\otimes_n (\nu_U\otimes_{i} \mu_U)}
    &\text{\scriptsize {[l]}}
    &
    &&&&\bullet
    \\
    &&&\text{\scriptsize {[j]}}
    &&\text{\scriptsize {[k]}}
    &&&\text{\scriptsize {[m]}}
    \\
    &&&&&\bullet
    \ar[rr]|{(1\otimes_{i} {\bcal J})\otimes_n (1\otimes_{i} 1)}
    &&\bullet
    \ar[dr]|{{\bcal M}\otimes_n {\bcal M}}
    \\
    &&&&\bullet
    \ar[ur]^{=}
    \ar[rrrr]|{1\otimes_n {\bcal M}}
    &&&&\bullet
    \ar[uuurrr]|{{\bcal {\bcal M}}}
    }
    $$
    \end{sideways}
  }
    	                    \end{center}
    	                    
  The arrows marked with an ``='' all occur as copies of $I$ are tensored to the object at the arrow's source. Therefore
    the quadrilateral regions [a],[e],[f],[h],[i] and [j] all commute trivially. The uppermost and lowermost quadrilaterals
    commute by the property of composing with units in an enriched $n$--category. The two triangles commute by the external
    unit condition for iterated monoidal categories. Regions [b] and [k] commute by respect of units by enriched 
$n$--functors.
    Regions [c] and [l] commute by naturality of $\eta.$ Region [g] commutes by the axiom of Definition~\ref{kcell}
   for $\nu$ and $\mu.$ Regions [d] and [m] commute by the ${\cal V}$--$n$--functoriality of ${\bcal M}.$ Note that the 
   instances of ${\bcal M}$ in these latter regions while written identically, actually have domains of differing categorical dimension.

The next case is $m=1$. The commuting diagram which shows that $\beta \circ \alpha$ obeys the axiom is as follows:
     
      \noindent
      \begin{center}
	          \resizebox{6in}{!}{ 
	      
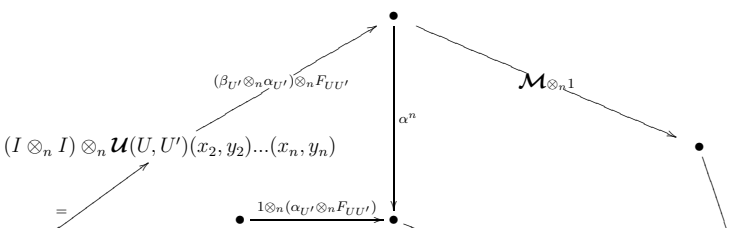
    \begin{sideways}\nopagebreak    
        $$
        \xymatrix
        @R=27pt@C=-64pt
        {
        &&&&\text{ }\text{ }\text{ }\text{ }\text{ }\text{ }\text{ }\text{ }\text{ }\text{ }\text{ }\text{ }\text{ }\text{ }\text{ }\text{ }\text{ }\text{ }\text{ }\text{ }\text{ }\text{ }
        \text{ }\text{ }\text{ }\text{ }\text{ }\text{ }\text{ }\text{ }\text{ }\text{ }\text{ }\text{ }\text{ }\bullet\text{ }\text{ }\text{ }\text{ }\text{ }\text{ }\text{ }\text{ }\text{ }\text{ }\text{ }
        \text{ }\text{ }\text{ }\text{ }\text{ }\text{ }\text{ }\text{ }\text{ }\text{ }\text{ }\text{ }\text{ }\text{ }\text{ }\text{ }\text{ }\text{ }\text{ }\text{ }\text{ }\text{ }\text{ }\text{ }
        \ar[ddrr]|{{\bcal M}\otimes_n 1}
        \ar[ddd]^{\alpha^n}
        \\\\
        &&(I \otimes_n I) \otimes_n {\bcal U}(U,U')(x_2,y_2)...(x_{n},y_{n})
        \ar[uurr]|{(\beta_{U'} \otimes_n \alpha_{U'})\otimes_n F_{UU'}}
        &&&&
        \text{ }\text{ }\text{ }\text{ }\text{ }\text{ }\text{ }\text{ }\text{ }\bullet\text{ }\text{ }\text{ }\text{ }\text{ }\text{ }\text{ }\text{ }\text{ }
        \ar[ddddr]^{{\bcal M}}
        \\
        &&&
        \bullet 
        \ar[r]^{1 \otimes_n (\alpha_{U'}\otimes_n F_{UU'})}
        &
        \bullet 
        \ar[dr]|{1\otimes_n {\bcal M}}
        \\
        &I \otimes_n {\bcal U}(U,U')(x_2,y_2)...(x_{n},y_{n})
        \ar[uur]^{=}
        \ar[dr]|{\beta_{U'} \otimes_n 1}
        &&&&
        \text{ }\text{ }\text{ }\text{ }\text{ }\text{ }\bullet\text{ }\text{ }\text{ }\text{ }\text{ }\text{ }
        \ar[ddrr]|{{\bcal M}}
        \ar@{=}[dd]
        \\
        &&\bullet
        \ar[uur]_{=}
        \ar[ddr]^{=}
        &&&
        \\
        &&&&&
        \text{ }\text{ }\text{ }\text{ }\text{ }\text{ }\bullet\text{ }\text{ }\text{ }\text{ }\text{ }\text{ }
        \ar[ddrr]|{{\bcal M}}
        &&
        \text{ }\text{ }\text{ }\text{ }\text{ }\bullet\text{ }\text{ }\text{ }\text{ }
        \ar@{=}[dd]
        \\
        &&&\bullet 
        \ar[r]^{1\otimes_n (G_{UU'}\otimes_n \alpha_U)}
        &\bullet 
        \ar[ur]|{1\otimes_n {\bcal M}}
        \\
        {\bcal U}(U,U')(x_2,y_2)...(x_{n},y_{n})
        \ar[uuuur]^{=}
        \ar[ddddr]_{=}
        &&&&&&&
        \text{ }\text{ }\text{ }\text{ }\text{ }\bullet\text{ }\text{ }\text{ }\text{ }
        \\
        &&&\bullet 
        \ar[r]^{(\beta_{U'}\otimes_n G_{UU'})\otimes_n 1}
        &\bullet
        \ar[dr]|{{\bcal M}\otimes_n 1}
        \ar[uu]_{\alpha^n}
        \\
        &&&&&
        \text{ }\text{ }\text{ }\text{ }\text{ }\text{ }\bullet\text{ }\text{ }\text{ }\text{ }\text{ }\text{ }
        \ar[uurr]|{{\bcal M}}
        \ar@{=}[dd]
        &&
        \text{ }\text{ }\text{ }\text{ }\text{ }\bullet\text{ }\text{ }\text{ }\text{ }
        \ar@{=}[uu]
        \\
        &&\bullet
        \ar[uur]_{=}
        \ar[ddr]^{=}
        &&&
        \\
        &{\bcal U}(U,U')(x_2,y_2)...(x_{n},y_{n}) \otimes_n I
        \ar[ur]|{1\otimes_n \alpha_U}
        \ar[ddr]_{=}
        &&&&
        \text{ }\text{ }\text{ }\text{ }\text{ }\bullet\text{ }\text{ }\text{ }\text{ }\text{ }\text{ }\text{ }
        \ar[uurr]|{{\bcal M}}
        \\
        &&&\bullet
        \ar[r]^{(H_{UU'}\otimes_n \beta_{U})\otimes_n 1}
        &\bullet
        \ar[ur]|{{\bcal M}\otimes_n 1}
        \ar[ddd]^{\alpha^n}
        \\
        &&{\bcal U}(U,U')(x_2,y_2)...(x_{n},y_{n})\otimes_n (I \otimes_n I)
        \ar[ddrr]|{H_{UU'}\otimes_n (\beta_{U}\otimes_n \alpha_U)}
        &&&&\text{ }\text{ }\text{ }\text{ }\text{ }\text{ }\text{ }\text{ }\text{ }\text{ }\text{ }\text{ }\text{ }
        \text{ }\text{ }\text{ }\text{ }\text{ }\text{ }\bullet\text{ }\text{ }\text{ }\text{ }\text{ }\text{ }\text{ }\text{ }\text{ }
        \text{ }\text{ }\text{ }\text{ }\text{ }\text{ }\text{ }\text{ }\text{ }\text{ }\text{ }\text{ }\text{ }\text{ }\text{ }
        \ar[uuuur]_{{\bcal M}}
        \\\\
        &&&&\text{ }\text{ }\text{ }\text{ }\text{ }\text{ }\text{ }\text{ }\text{ }\text{ }\text{ }\text{ }\text{ }\text{ }\text{ }\text{ }
        \text{ }\text{ }\text{ }\text{ }\text{ }\text{ }\text{ }\text{ }\text{ }\text{ }\bullet\text{ }\text{ }\text{ }\text{ }\text{ }\text{ }\text{ }\text{ }
        \text{ }\text{ }\text{ }\text{ }\text{ }\text{ }\text{ }\text{ }\text{ }\text{ }\text{ }\text{ }\text{ }\text{ }\text{ }\text{ }\text{ }\text{ }
        \ar[uurr]|{1\otimes_n {\bcal M}}
        }
        $$\nopagebreak
         \end{sideways}
         }
		                    \end{center}
		
The arrows marked with an ``='' all occur as copies of $I$ are tensored to the object at the arrow's source. 
         The 3 leftmost regions commute by the naturality of $\alpha^n.$ The 2 embedded central ``hexagons'' commute by the
         axiom in Definition~\ref{kcell} for $\beta$ and $\alpha.$ The three pentagons on the right
         are copies of the pentagon axiom for the composition ${\bcal M}.$               

Finally, the last case we examine is $m=0$. To reduce it to the previous case we apply the principle of whiskering
of $(n+1)$--cells with ${\cal V}$--$n$--functors. Once we have done so there are two obvious ways of describing the
composition of two $(n+1)$--cells along a common ${\cal V}$--$n$--category in terms of composition along a common
${\cal V}$--$n$--functor. The two ways are equivalent since everything is taking place in the larger $(n+1)$--category
structure. Now there are two subcases corresponding to
the left and the right-hand whiskering. Closure under whiskering onto the right is shown by:
\noindent
        		          \begin{center}
	          \resizebox{6.55in}{!}
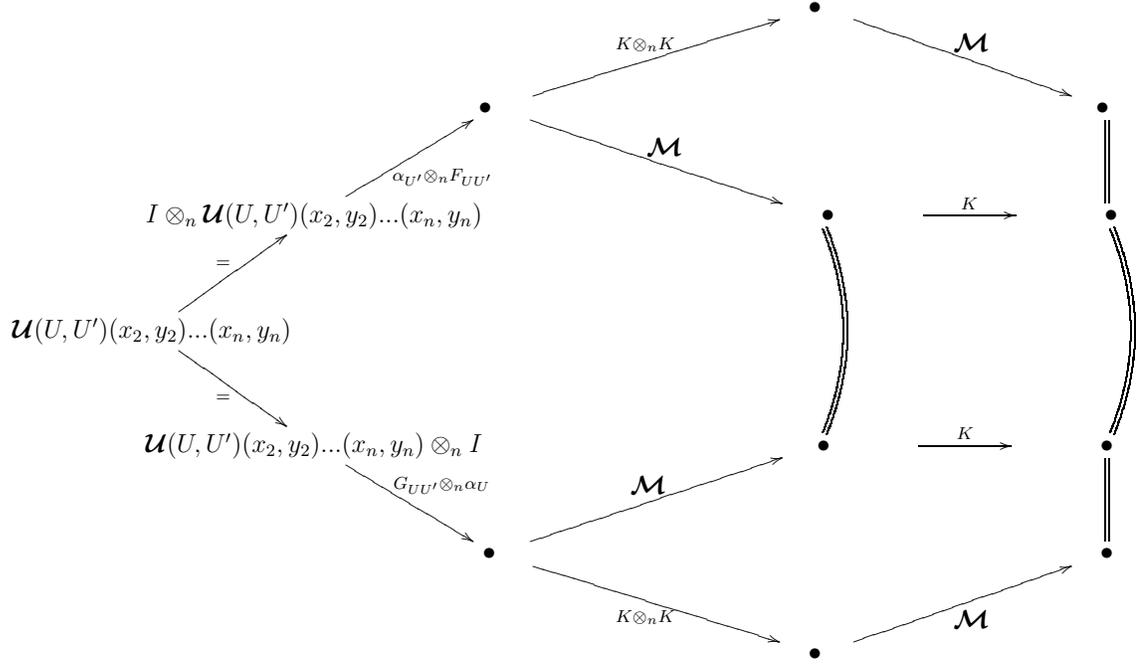
{ 
  $$
             \xymatrix@C=-73pt@R=35pt{
             &&&\text{ }\text{ }\text{ }\text{ }\text{ }\text{ }\text{ }\text{ }\text{ }\text{ }\text{ }\text{ }\text{ }\text{ }\text{ }\text{ }\text{ }\text{ }\text{ }\text{ }\text{ }\text{ }\text{ }
             \text{ }\text{ }\text{ }\text{ }\text{ }\text{ }\text{ }\text{ }\text{ }\text{ }\text{ }\text{ }\bullet\text{ }\text{ }\text{ }\text{ }\text{ }\text{ }\text{ }\text{ }\text{ }\text{ }\text{ }\text{ }
             \text{ }\text{ }\text{ }\text{ }\text{ }\text{ }\text{ }\text{ }\text{ }\text{ }\text{ }\text{ }\text{ }\text{ }\text{ }\text{ }\text{ }\text{ }\text{ }\text{ }\text{ }\text{ }\text{ }\text{ }
             \ar[dr]^-{{\bcal M}}
             \\
             &
             &\text{ }\text{ }\text{ }\text{ }\text{ }\text{ }\text{ }\text{ }\text{ }\text{ }
             \text{ }\text{ }\text{ }\text{ }\text{ }\text{ }\bullet\text{ }\text{ }\text{ }\text{ }\text{ }\text{ }
             \text{ }\text{ }\text{ }\text{ }\text{ }\text{ }\text{ }\text{ }\text{ }\text{ }\text{ }\text{ }
             \ar[rd]^-{{\bcal M}}
             \ar[ru]^{K\otimes_n K}
             &&\text{ }\text{ }\text{ }\text{ }\text{ }\text{ }\text{ }
             \text{ }\text{ }\text{ }\text{ }\text{ }\bullet\text{ }\text{ }\text{ }\text{ }
             \text{ }\text{ }\text{ }\text{ }\text{ }\text{ }\text{ }\text{ }\text{ }
             \ar@{=}[d]
           \\
           &
           I \otimes_n {\bcal U}(U,U')(x_2,y_2)...(x_{n},y_{n})
             \ar[ru]_{\text{ }\alpha_{U'} \otimes_n F_{UU'}}
           &&\text{ }\text{ }\text{ }\text{ }\text{ }\text{ }\text{ }\text{ }
           \text{ }\text{ }\text{ }\bullet\text{ }\text{ }\text{ }
           \text{ }\text{ }\text{ }\text{ }\text{ }\text{ }
           \ar@{=}@/^1pc/[dd]
           \ar[r]^-{K}
           &\text{ }\text{ }\text{ }\text{ }\text{ }\text{ }
           \text{ }\text{ }\text{ }\bullet\text{ }\text{ }\text{ }
           \text{ }\text{ }\text{ }\text{ }\text{ }
           \ar@{=}@/^1pc/[dd]
           \\
           {\bcal U}(U,U')(x_2,y_2)...(x_{n},y_{n})
             \ar[ru]^{=}
           \ar[rd]_{=}
           \\
           &
           {\bcal U}(U,U')(x_2,y_2)...(x_{n},y_{n}) \otimes_n I
             \ar[rd]^{\text{ }G_{UU'} \otimes_n \alpha_{U}}
           &&\text{ }\text{ }\text{ }\text{ }\text{ }\text{ }\text{ }
           \text{ }\text{ }\text{ }\bullet\text{ }\text{ }\text{ }
           \text{ }\text{ }\text{ }\text{ }\text{ }\text{ }
           \ar[r]^-{K}
           &\text{ }\text{ }\text{ }\text{ }\text{ }\text{ }
           \text{ }\text{ }\text{ }\bullet\text{ }\text{ }\text{ }
           \text{ }\text{ }\text{ }\text{ }\text{ }\text{ }
           \\
             &
             &\text{ }\text{ }\text{ }\text{ }\text{ }\text{ }\text{ }\text{ }\text{ }
             \text{ }\text{ }\text{ }\text{ }\bullet\text{ }\text{ }\text{ }\text{ }\text{ }
             \text{ }\text{ }\text{ }\text{ }\text{ }\text{ }\text{ }\text{ }\text{ }
             \ar[ru]^-{{\bcal M}}
             \ar[rd]_{K\otimes_n K}
             &&\text{ }\text{ }\text{ }\text{ }\text{ }\text{ }
             \text{ }\text{ }\text{ }\bullet\text{ }\text{ }\text{ }
             \text{ }\text{ }\text{ }\text{ }\text{ }\text{ }
             \ar@{=}[u]
             \\
             &&&\text{ }\text{ }\text{ }\text{ }\text{ }\text{ }\text{ }\text{ }\text{ }\text{ }
             \text{ }\text{ }\text{ }\text{ }\text{ }\bullet\text{ }\text{ }\text{ }\text{ }\text{ }\text{ }\text{ }
             \text{ }\text{ }\text{ }\text{ }\text{ }\text{ }\text{ }\text{ }\text{ }
             \ar[ur]_-{{\bcal M}}
             }
           $$
           }
  \end{center}
The preceding diagram's commutativity
  relies on the fact that $\alpha$ obeys Definition~\ref{kcell} and on the ${\cal V}$--n--functoriality of $K.$
 
Closure under whiskering onto the left is shown by the following
commutative diagram:
  
\noindent
      		          \begin{center}
	          \resizebox{6.55in}{!}
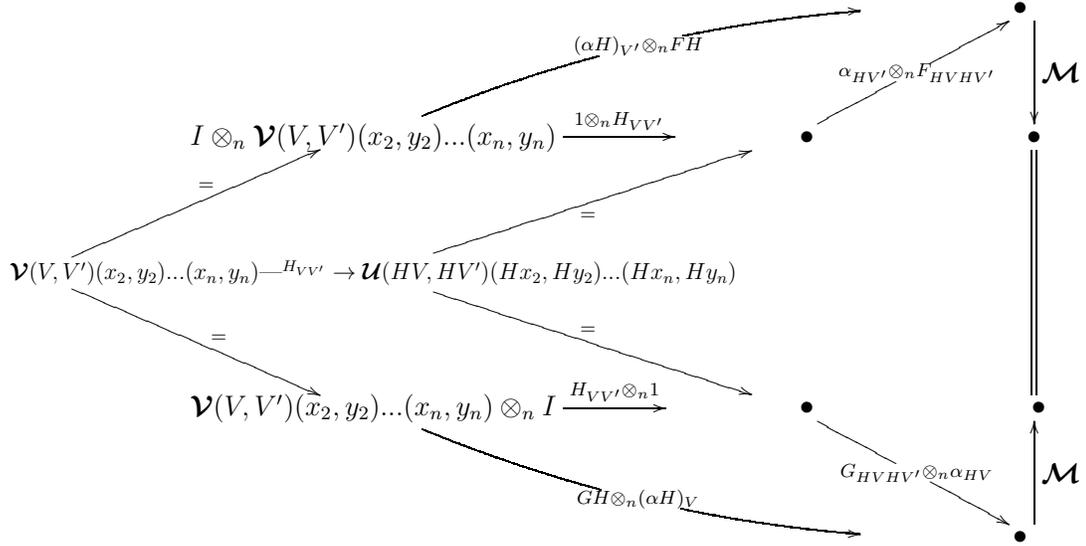
{ 
  $$
             \xymatrix@C=-33pt@R=40pt{
             &
             &&\text{ }\text{ }\text{ }\text{ }\text{ }\text{ }\text{ }\text{ }\text{ }\text{ }\text{ }\text{ }
             \text{ }\text{ }\text{ }\bullet\text{ }\text{ }\text{ }\text{ }\text{ }\text{ }
             \text{ }\text{ }\text{ }\text{ }\text{ }\text{ }\text{ }\text{ }\text{ }\text{ }\text{ }\text{ }
             \ar[d]^-{{\bcal M}}
           \\
           &
           I\otimes_n {\bcal V}(V,V')(x_2,y_2)...(x_{n},y_{n})
           \ar[r]^-{1\otimes_n H_{VV'}}
           \ar@/^1.4pc/[rru]|-{(\alpha H)_{V'} \otimes_n FH}
           &
           \text{ }\text{ }\text{ }\text{ }\text{ }\text{ }\text{ }\text{ }\text{ }
           \text{ }\text{ }\text{ }\bullet\text{ }\text{ }\text{ }
           \text{ }\text{ }\text{ }\text{ }\text{ }\text{ }
             \ar[ru]|-{\text{ }\alpha_{HV'} \otimes_n F_{HVHV'}}
           &\text{ }\text{ }\text{ }\text{ }\text{ }\text{ }\text{ }\text{ }\text{ }\text{ }\text{ }\text{ }
           \text{ }\text{ }\text{ }\text{ }\text{ }\text{ }\bullet\text{ }\text{ }\text{ }\text{ }\text{ }\text{ }
           \text{ }\text{ }\text{ }\text{ }\text{ }\text{ }\text{ }\text{ }\text{ }\text{ }\text{ }\text{ }
           \ar@{=}[dd]
           &
           \\
           \text{ }\text{ }\text{ }
           \ar@<4pt>[ru]^{=}
           \ar@<-4pt>[rd]^{=}
           &\text{\resizebox{4.25in}{!}{${\bcal V}(V,V')(x_2,y_2)...(x_{n},y_{n})\text{---}{}^{H_{VV'}}\to{\bcal U}(HV,HV')(Hx_2,Hy_2)...(Hx_n,Hy_n)$}}
           \ar[ru]_{=}
           \ar[rd]^{=}
           &
           \\
           &
           {\bcal V}(V,V')(x_2,y_2)...(x_{n},y_{n})\otimes_n I
  	 \ar[r]^-{H_{VV'}\otimes_n 1}
           \ar@/_1.4pc/[rrd]|-{GH \otimes_n (\alpha H)_V}
           &
           \text{ }\text{ }\text{ }\text{ }\text{ }\text{ }\text{ }\text{ }
           \text{ }\text{ }\text{ }\text{ }\text{ }\bullet\text{ }\text{ }\text{ }\text{ }
           \text{ }\text{ }\text{ }\text{ }\text{ }\text{ }
             \ar[rd]|-{\text{ }G_{HVHV'} \otimes_n \alpha_{HV}}
           &\text{ }\text{ }\text{ }\text{ }\text{ }\text{ }\text{ }\text{ }\text{ }\text{ }\text{ }\text{ }\text{ }
           \text{ }\text{ }\text{ }\text{ }\text{ }\text{ }\bullet\text{ }\text{ }\text{ }\text{ }\text{ }\text{ }
           \text{ }\text{ }\text{ }\text{ }\text{ }\text{ }\text{ }\text{ }\text{ }\text{ }\text{ }\text{ }
           &
           \\
             &&
             &\text{ }\text{ }\text{ }\text{ }\text{ }\text{ }\text{ }\text{ }\text{ }
             \text{ }\text{ }\text{ }\text{ }\text{ }\text{ }\bullet\text{ }\text{ }\text{ }\text{ }\text{ }\text{ }
             \text{ }\text{ }\text{ }\text{ }\text{ }\text{ }\text{ }\text{ }\text{ }\text{ }\text{ }\text{ }
             \ar[u]_-{{\bcal M}}
             }
           $$
           }
  \end{center}
  The exterior of this diagram commutes simply because $\alpha$ obeys the axiom of Definition~\ref{kcell}.
                    

There are specific instances of all the above constructions and especially of these last few commuting diagrams found in 
\cite{forcey2}. 

In conclusion, it should be clear that the enriched $n$--categories we have constructed generalize the 
notion of strict $n$--category. When a complete definition of weak $n$--category
is in place there should be an analogous generalization to weak enrichment. Another direction to go in is that of
${\cal V}$--Mod, the category of ${\cal V}$--categories with ${\cal V}$--modules as morphisms.
${\cal V}$--Mod should also be investigated for the case of ${\cal V}$ $k$--fold monoidal. The same is true of 
${\cal V}$--Act, the category of categories with a ${\cal V}$ action. In any case, further study should attempt to 
elucidate the relationship of the nerves of the $n$--categories in question. For instance, we would like to know
the relationship between $\Omega Nerve$(${\cal V}$--Cat) and $Nerve({\cal V}).$  This would even be quite interesting
in the case of symmetric ${\cal V}$ where the nerve is an infinite loop space. It would be nice to know if there are
symmetric monoidal categories whose nerves exhibit periodicity under the vertically iterated enrichment functor.
    

        
    }
    
\end{document}